\documentclass[final,onefignum,onetabnum]{siamonline220329}

\usepackage{mmap}

\usepackage{braket,amsfonts}

\usepackage{array}

\usepackage[caption=false]{subfig}
\usepackage{pgfplots}

\newsiamthm{claim}{Claim}
\newsiamremark{remark}{Remark}
\newsiamremark{hypothesis}{Hypothesis}
\crefname{hypothesis}{Hypothesis}{Hypotheses}

\usepackage{algorithmic}

\usepackage{graphicx,epstopdf}

\Crefname{ALC@unique}{Line}{Lines}

\usepackage{amsopn}

\newcommand{\Path}{}
\usepackage{enumerate}
\usepackage{tikz-cd}
\usepackage{natbib}
\newcommand{\figs}{}

\newcommand{\reconstruct}{\mathcal{T}}

\newcommand{\cocyc}{\mathcal{G}}

\newcommand{\W}{\mathfrak{S}}

\newcommand{\Csens}{ C_{\phi,\Phi} }
\DeclareMathOperator{\retract}{ret}

\DeclareMathOperator{\AutCor}{AutCorr}
\newcommand{\calH}{\mathcal{H}}
\newcommand{\calV}{\mathcal{V}}
\newcommand{\calD}{\mathcal{D}}
\newcommand{\calU}{\mathcal{U}}
\newcommand{\calX}{\mathcal{X}}
\newcommand{\calW}{\mathcal{W}}
\newcommand{\calP}{\mathcal{P}}
\DeclareMathOperator{\clos}{clos}
\definecolor{PineGreen}{rgb}{0.0, 0.3, 0.2}
\definecolor{Itranga}{rgb}{0.99, 0.27, 0.0}
\newcommand{\itranga}{\textcolor{Itranga}}
\definecolor{Mete}{RGB}{92, 75, 62}

\definecolor{Turquoise}{RGB}{22, 219, 219}
\definecolor{Akashi}{RGB}{115, 228, 245}

\newcommand{\blue}{\textcolor{blue}}

\definecolor{ChhaiA}{rgb}{0.9, 0.9, 0.9}
\definecolor{ChhaiB}{rgb}{0.8, 0.8, 0.8}
\definecolor{ChhaiC}{rgb}{0.7, 0.7, 0.7}
\definecolor{ChhaiD}{rgb}{0.6, 0.6, 0.6}
\definecolor{ChhaiE}{rgb}{0.5, 0.5, 0.5}
\definecolor{ChhaiF}{rgb}{0.4, 0.4, 0.4}
\definecolor{ChhaiG}{rgb}{0.3, 0.3, 0.3}
\definecolor{ChhaiH}{rgb}{0.2, 0.2, 0.2}
\definecolor{SheolaA}{rgb}{0.69, 0.9, 0.7}
\definecolor{SheolaB}{rgb}{0.38, 0.47, 0.4}
\newtheorem{Assumption}{Assumption}
\usepackage{tabularx} \newcolumntype{L}{>{\raggedright\arraybackslash}X}
\newcommand{\real}{\mathbb{R}}

\newcommand{\integer}{\mathbb{Z}}
\newcommand{\num}{\mathbb{N}}
\newcommand{\Disc}{\mathcal{D}}
\newcommand{\bigO}[1]{O \left( #1 \right)}

\newcommand{\Linear}{\mathcal{L}}
\DeclareMathOperator{\proj}{proj}
\DeclareMathOperator{\ran}{ran}

\DeclareMathOperator{\Id}{Id}
\DeclareMathOperator{\spn}{span}
\newcommand{\imply}{\quad \Rightarrow \quad}

\newcommand{\bracketBig}[1]{\left\langle #1 \right\rangle}
\newcommand{\norm}[1]{\left\| #1 \right\|}

\newcommand{\SetDef}[2]{ \left\{ #1 \;:\; #2  \right\} }
\newcommand{\paran}[1]{\left( #1 \right)}
\newcommand{\braces}[1]{\left\{ #1 \right\}}

\usepackage{xspace}
\usepackage{bold-extra}
\usepackage[most]{tcolorbox}

\colorlet{texcscolor}{blue!50!black}
\colorlet{texemcolor}{red!70!black}
\colorlet{texpreamble}{red!70!black}
\colorlet{codebackground}{black!25!white!25}
\usepackage{tikz-cd}
\usepackage{\Path tikzcd2}

\lstdefinestyle{siamlatex}{%
	style=tcblatex,
	texcsstyle=*\color{texcscolor},
	texcsstyle=[2]\color{texemcolor},
	keywordstyle=[2]\color{texemcolor},
	moretexcs={cref,Cref,maketitle,mathcal,text,headers,email,url},
}

\tcbset{%
	colframe=black!75!white!75,
	coltitle=white,
	colback=codebackground, 
	colbacklower=white, 
	fonttitle=\bfseries,
	arc=0pt,outer arc=0pt,
	top=1pt,bottom=1pt,left=1mm,right=1mm,middle=1mm,boxsep=1mm,
	leftrule=0.3mm,rightrule=0.3mm,toprule=0.3mm,bottomrule=0.3mm,
	listing options={style=siamlatex}
}

\newtcblisting[use counter=example]{example}[2][]{%
	title={Example~\thetcbcounter: #2},#1}

\newtcbinputlisting[use counter=example]{\examplefile}[3][]{%
	title={Example~\thetcbcounter: #2},listing file={#3},#1}

\DeclareTotalTCBox{\code}{ v O{} }
{ 
	fontupper=\ttfamily\color{black},
	nobeforeafter,
	tcbox raise base,
	colback=codebackground,colframe=white,
	top=0pt,bottom=0pt,left=0mm,right=0mm,
	leftrule=0pt,rightrule=0pt,toprule=0mm,bottomrule=0mm,
	boxsep=0.5mm,
	#2}{#1}

\patchcmd\newpage{\vfil}{}{}{}
\title{Limits of Learning Dynamical Systems
	\thanks{\funding{Both authors acknowledge support by NSF grant DMS-204839.  Tyrus Berry also acknowledges support by NSF grant DMS-2006808.}}}

\author{Tyrus Berry \thanks{Department of Mathematical Sciences, George Mason University (\email{tyrus.berry@gmail.com}).} \orcid{https://orcid.org/0000-0002-0280-7508}
	\and Suddhasattwa Das \thanks{Department of Mathematics and Statistics, Texas Tech University (\email{iamsuddhasattwa@gmail.com}).} \orcid{https://orcid.org/0000-0003-2085-4743}}

\headers{Limits of learning Dynamical Systems}{Tyrus Berry and Suddhasattwa Das}

\begin{document}
	\slugger{siads}{XXXX}{0}{0}{000--000}
	\maketitle
	
	\begin{abstract} A dynamical system is a transformation of a phase space, and the transformation law is the primary means of defining as well as identifying the dynamical system. It is the object of focus of many learning techniques. Yet there are many secondary aspects of dynamical systems - invariant sets, the Koopman operator, and Markov approximations, which provide alternative objectives for learning techniques. Crucially, while many learning methods are focused on the transformation law, we find that forecast performance can depend on how well these other aspects of the dynamics are approximated.  These different facets of a dynamical system correspond to objects in completely different spaces - namely interpolation spaces, compact Hausdorff sets, unitary operators and Markov operators respectively. Thus learning techniques targeting any of these four facets perform different kinds of approximations. We examine whether an approximation of any one of these aspects of the dynamics could lead to an approximation of another facet. Many connections and obstructions are brought to light in this analysis. Special focus is put on methods of learning of the primary feature - the dynamics law itself. The main question considered is the connection of learning this law with reconstructing the Koopman operator and the invariant set. The answers are tied to the ergodic and topological properties of the dynamics, and reveal how these properties determine the limits of forecasting techniques.
	\end{abstract}
	
	
	
	\section{Introduction} \label{sec:intro}
	
	Many real world phenomenon are modeled as dynamical systems, and these dynamical systems models often need to be discovered from data. There are many disciplines such as climate sciences \cite{MessoriFaranda2021paleo, ghil2017wind, SlawinskaGiannakis2017}, traffic dynamics \cite{DasEtAl2023traffic, DasMustAgar2023_qpd} and epidemiology \cite{piazzola2021note, MustaveeEtAl_covid_2021} where data-driven discovery is an important step towards understanding the system. With the growth of computational power as well as data, many new techniques and paradigms of reconstructing dynamical systems have emerged. Indeed, while some data types (eg. images and language) seem to be converging on a preferred architecture (eg. convolutional networks and transformers respectively); on the contrary, there appears to be no universal architecture for learning dynamical data.  While neural network methods such as reservoir computers and long short-term memory (LSTM) architectures have had considerable success on many problems, their architectures are fundamentally different and other learning techniques such as dynamic mode decomposition, embedology based methods, and stochastic approximations are still competitive on many problems. The reason these diverse techniques have persisted is that they each capture different aspects of the dynamics. To develop effective machine learning methods for dynamical data, one must understand the role played by various properties of the underlying dynamics.
	
\begin{figure}\center
		\begin{tikzpicture}[scale=0.6, transform shape]
			\input{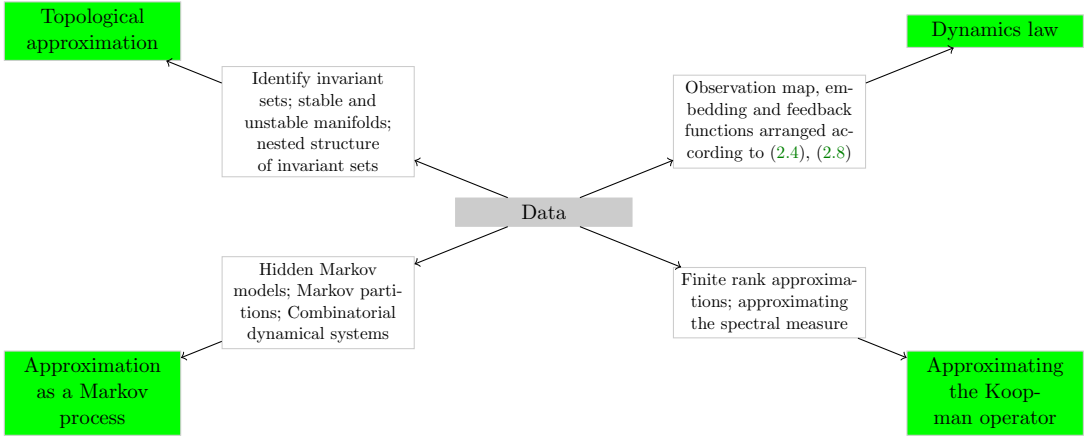}
		\end{tikzpicture}
		\caption{Many features of a dynamical system. If a dynamical system generates data, there are many different aspects of the dynamics that one can try to recover from the data. These aspects are shown in green, and some common techniques in white. }
		\label{fig:outline1}
	\end{figure}

	A dynamical system can be classically defined as a group of self-maps on a set, usually continuous maps on a topological space. There are also many alternative means of describing the dynamics which are more relevant in other contexts. Figure \ref{fig:outline1} presents an overview of these aspects of a dynamical system. As a result, there have been many techniques targeting these different aspects, each having its own set of advantages and disadvantages. Our goal is to classify these techniques and also look in depth at their inter-relations. These inter-relations provide an organized overview of the existing state of knowledge in the development of learning techniques for dynamical systems. They also provide an overview of several important conceptual gaps in this field. 
	
	Our focus is on discrete-time dynamics. Formally, we assume :
	
	\begin{Assumption} \label{A:f}
		There is a $C^1$ dynamical system $f: \tilde\Omega \to \tilde\Omega$ on an $m$-dimensional $C^1$ manifold $\tilde\Omega$, with an ergodic measure $\mu$ with compact support $\Omega$.
	\end{Assumption}
	
	This assumption is general enough to capture the essence of many systems, as well as modeling a dynamical source that generates a stream of data. The assumption of an ergodic invariant measure with compact support is satisfied in any system with bounded trajectories. We next describe how data in the form of timeseries may be generated from such systems.
	
	\paragraph{Data-driven discovery} In the various techniques that we shall describe, the ingredients in Assumption \ref{A:f} are usually unknown. Their existence is a part of the null hypothesis. Our only access to these unknown ingredients will be through a measurement (also known as an observation). 
	
	\begin{Assumption} \label{A:data}
		There is a $C^1$ function $\phi : \tilde\Omega \to \real^d$, serving as a measurement of the system.
	\end{Assumption}
	
	Assumptions \ref{A:f} and \ref{A:data} provide the basis for a data-driven study of dynamical systems. Given any initial state $\omega_0 \in \Omega$, one obtains a time series
	\begin{equation} \label{eqn:def:hatomega}
		\braces{ \hat{\omega}_n }_{n\in\num_0} , \quad \hat{\omega}_n := \phi \paran{ f^n \omega_0 } .
	\end{equation}
	Conversely, any timeseries $\braces{ \hat{\omega}_n }_{n\in\num_0}$ can be assumed to be originating from a system satisfying Assumptions \ref{A:f} and \ref{A:data}. One can then try to infer $f$ and $\phi$ from the timeseries. This analysis of the timeseries becomes a classic problem in learning theory or functional approximation. 
	
	This task of learning the dynamics law is presented as one of the four branches in Figure \ref{fig:outline1}. It is one of the many aspects of the broader field of \emph{data-driven discovery} of dynamical systems. The dynamics law $f$ is one of the many possible ways to characterize the underlying dynamics, although it serves as our primary description of the system. Each aspect of the dynamics presented in Figure \ref{fig:outline1} provides a separate objective for numerical methods designed for dynamical systems, and for data-driven learning methods in particular. A detailed mathematical description of this task is presented in Section \ref{sec:learn}. We now briefly describe the three other branches of Figure \ref{fig:outline1}.
	
	\paragraph{The invariant set} A subset $X\subseteq \tilde\Omega$ is said to be invariant if $f(X) \subseteq X$. Assumption \ref{A:f} specifies such an invariant set $\Omega$. Any discussion of a dynamical system is usually focused on one or more invariant sets. An invariant set could be a fixed point, a periodic orbit, a Hamiltonian torus or a strange attractor as seen in chaotic systems. An invariant set is a subset of the total phase space which forms its own dynamical system under the same transition law $f$. It presents its own local topological and measure theoretic features which may not be valid globally.
	If an orbit starts on an invariant set $X$, it remains in $X$ for all times. Thus the timeseries it generates provides information limited to $X$. If the measurement $\phi$ from Assumption \ref{A:data} is injective, one obtains an embedding of $\Omega$ in $\real^d$. The data points in the timeseries thus provide a sampling of $\Omega$. More importantly, they do not sample regions in the phase space lying outside $\Omega$. The inferences we draw from numerical methods will be limited to the subset being sampled. 
	
	Given any point $z\in \tilde\Omega$, its $\omega$-limit set is the collection of accumulation points of the sequence $\braces{ f^n z }_{n=0}^{\infty}$. It is thus the collection of all those points which are approached arbitrarily closely infinitely many times by the trajectory of $z$. An $\omega$-limit set can be shown to be an invariant closed set itself. The \emph{basin} of an invariant set $X$ is the set of points $z$ whose $\omega$-limit set is contained within $X$. In other words, it is the collection of points which converge towards $X$ as the dynamics unfolds. Given a point $a\in \tilde\Omega$ let $\delta_{a}$ denote the Dirac-delta measure, which is the Borel probability measure whose support is the singleton set $\{a\}$. Note that as the point $z$ travels along its trajectory one obtains a sequence of measures $\frac{1}{N} \sum_{n=0}^{N-1} \delta_{f^n z}$. These are discrete measures called \emph{empirical measures}, as they are sample based. A point $z$ is said to be in the basin of the invariant measure $\mu$ if its empirical measures converge weakly to $\mu$. We will call the set $X$ or an invariant measure $\mu$ to be \emph{visible} if their basins have non-zero volume. Thus being visible means that they may be approximated by observations made from physical experiments. We shall call $X$ or $\mu$ \emph{stable} if their basins contain a neighborhood of $X$ and the support of $\mu$ respectively. In these cases $X$ and $\mu$ are also called an \emph{attractor} and a \emph{physical measure} respectively. This means that given any point $z$ close enough to either of these sets, the trajectory of $z$ provides a means of approximating $X$ or $\mu$. Visibility and stability are important properties that allow phenomena and subsets of the phase-space to be accessible to numerical experiments and observations. These will become important in subsequent sections.
	
	If an invariant set $X$ is stable it is also called an \emph{attractor}. Its basin contains an open neighborhood $\calU$ of $X$. Then we have
	\begin{equation} \label{eqn:X_U}
		X = \cap_{n\in\num_0} \clos f^n(\calU) .
	\end{equation}
	Alternatively if there is an open region $\calU$ in $\tilde\Omega$ such that $\clos f(\calU) \subset \calU$, then one can apply \eqref{eqn:X_U} to obtain a closed invariant attracting region $X$. Equation \eqref{eqn:X_U} re-interprets an attractor as the limiting set of an iterative procedure. This is more important from a computational perspective, as one usually does not have access to the attractor $X$ but to an invariant open domain $\calU$. The interpretation of an attractor as an asymptotic property also underlies its importance in understanding the nature of the dynamics.
	
	There have been many approaches to identifying an attractor or a general invariant set $\Omega$ as in Assumption \ref{A:f}. The samples collected from the set $\Omega$ can be fed into any topological data analysis technique \cite[e.g.]{BerrySauer2016, VaughnBerryAntil2019, BauerLesnick2020prsstnc, OtterEtAl2017roadmap} to obtain an estimate of $\Omega$. One can also take a structuralist approach to $\Omega$ and decompose it into chain-recurrent sets connected with stable and unstable manifolds \cite[e.g.]{RobbinSalamon1988, Mischaikow1999Conley}. Each chain recurrent set is invariant and allows a similar structural decomposition.
	
	\paragraph{Koopman operator} A general measurable function $\psi : \tilde\Omega \to \real$ is called an \emph{observable}. The Koopman operator $U$ \cite{Halmos1956, DasGiannakis_delay_2019, DasGiannakis_RKHS_2018} characterizes the dynamics not in terms of the trajectories on the phase space $\tilde\Omega$ but in terms of the evolution of observables. Any function $\psi : \tilde\Omega \to \real$ is transformed by $U$ into
	\begin{equation} \label{eqn:def:Koop}
		U\psi : \tilde\Omega \to \real, \; (U\psi)(\omega) = \psi \paran{ f \omega } .
	\end{equation}
	Thus $U$ is a linear operator on the space of functions. It acts on observables by composition with the map $f$. It has many invariant subspaces, such as $C^r(\tilde\Omega)$, $C(\Omega)$ or $L^2(\mu)$. The action and properties of $U$ can be studied separately on each of these function spaces. From an operator theoretic point of view, this most common choice of functional space is $L^2(\mu)$ -- the collection of square-integrable, complex-valued functions. This is a Hilbert space on which $U$ is unitary. The Koopman operator allows the use of various tools from functional analysis, otherwise unavailable in the classical approach to nonlinear dynamics.
	
	The Koopman operator provides an equation free description of the dynamics. It can be directly interpreted as the action of forecasting the future states of a measurement. In many application areas \cite[e.g.]{DasEtAl2023traffic, MustaveeEtAl_covid_2021, GiannakisDas_tracers_2019} it suffices to directly approximate the Koopman operator instead of the dynamics map $f$. The unitarity of $U$ is of great importance in the context of forecasting. Suppose there is a target function $\phi$ and an estimate $\hat\phi$ of $\phi$. These two functions would evolve under the dynamics as $U^t \phi$ and $U^t \hat\phi$. Due to unitarity we have
	\begin{equation} \label{eqn:untry:1}
		\norm{ U^t \phi - U^t \hat\phi }_{L^2(\mu)} = \norm{ \phi - \hat\phi }_{L^2(\mu)}, \quad \forall t\in \real.
	\end{equation}
	Equation \eqref{eqn:untry:1} indicates that any initial error in estimation does not grow in $L^2$-norm. Thus the Koopman operator is a stable predictor on $L^2(\mu)$.
	
	Of course, equation \eqref{eqn:untry:1} assumes an accurate representation of $U$ itself, but any technique for approximating $U$ must rely on finite matrices. Finite matrices or finite rank operators are often inadequate for approximating unitary operators on infinite dimensional Hilbert spaces. One way of tackling this problem is to look for an invariant subspace $\calV$ of $L^2(\mu)$ and simplify the task to characterizing the action of $U$ on $\calV$. Given any target function $\phi$ to forecast, if $\hat\phi$ denotes the orthogonal projection $\phi$ onto $\calV$, then \eqref{eqn:untry:1} says that the forecasting error will remain constant and equal to the norm of the residual of the projection. A common choice for $\calV$ is $\calD$, the closure of all eigenfunctions. There has been a number of techniques proposed \cite[e.g.]{DasGiannakis_delay_2019, DasGiannakis_RKHS_2018} which provide a convergent estimate of the discrete spectrum and functional space of $U$. A more general technique is presented in \cite{DGJ_compactV_2018, ValvaGiannakis2023cnstnt} for continuous time systems, which performs a spectral approximation of the generator $V$ of a continuous-time dynamical system. Approximation of the Koopman operator remains a challenging task, particularly in systems with different time-scales.
	
	\paragraph{Markov process}  A Markov process is a stochastic process $\SetDef{X_t}{t\geq 0}$ in which given a sequence of observation times $t_1 < t_2 < \ldots < t_n$, the current state at time $t_n$ depends only on the state at time $t_{n-1}$. More precisely we have
	\[ \mathbb{E} \paran{ X_{t_n} \,\rvert\, X_{t_{n-1}}, \ldots, X_{t_1} } = \mathbb{E} \paran{ X_{t_n} \,\rvert\, X_{t_{n-1}} } . \]
	Markov processes model almost all phenomena with a \emph{memory-less} property. These include the solutions of stochastic differential equations (SDE)s which are used to model phenomena such as stock market movements, as well as analytic studies such as error analysis of computer arithmetic. Markov processes are easy to simulate \cite[e.g.]{desharnais2002bisim, whitt1992asymp} and provide an alternative description of the trajectories of the dynamics. 
	
	Before attempting a description of a dynamical system as a Markov process, one needs to interpret it as a probability space. The invariant measure $\mu$ of the dynamics provides that interpretation. Then any function such as $\phi$ on $\tilde\Omega$ can be interpreted as a random variable. Moreover, any transform of $\phi$, such as a a time-average, or Koopman transform, is again a new variable. Given a sequence of observation times $t_1 < t_2 < \ldots < t_n$, one could ask whether the following equality holds :
	\[ \mathbb{E} \paran{ U^{t_n} \phi \,|\, U^{t_{n-1}} \phi, \ldots , U^{t_1} \phi   } \stackrel{?}{=} \mathbb{E} \paran{ U^{t_n} \phi \,|\, U^{t_{1}} \phi } . \]
	There is no definitive answer to this question. The question can be made easier by demanding that the process is a finite-state process. The states correspond to the cells of a finite partition $\calU = \SetDef{U_i}{ 1\leq i \leq m }$ of the invariant set. One requires that
	\[ \mathbb{E} \paran{ x_{t_n} \in U_{i_n} \,|\, x_{t_{n-1}} \in U_{i_{n-1}}, \ldots , x_{t_{1}} \in U_{i_{1}} } = \mathbb{E} \paran{ x_{t_n} \in U_{i_n} \,|\, x_{t_{n-1}} \in U_{i_{n-1}} } . \]
	Such a partition is known as a \emph{Markov partition}. The existence of Markov partitions has only been proved for special systems such as hyperbolic systems \cite[e.g.]{Bowen_Markov_1970}. A major improvement was showing that more general and abundant class of dynamics called \emph{non-uniformly hyperbolic} systems can be approximated by Markov partitions \cite{Sanchez_pressure_2017, LiaoEtAl2018}. There has been no general existence results known yet.
	
	Given any general partition $\calU$ as above, one obtains a stochastic process on the finite state space. There has been many algebraic characterization of such processes to be Markovian \cite[e.g.]{vidyasagar2005rlze, vanluyten2006matrix, rabiner1986intro, picci2005weak, Picci1978internal, Heller1965stoch, finesso2010approx}. One can avoid these issues and create an $m\times m$ \emph{transition matrix} $\mathbb{P}$ whose entries are 
	\begin{equation} \label{eqn:def:transition}
		\mathbb{P}_{i,j} := \mathbb{E} \paran{ f(x) \in U_{i} \,|\, x \in U_j } , \quad \forall 1\leq i,j \leq m.
	\end{equation}
	Such a matrix is column-stochastic and can be used to simulate a Markov process on its state space. This Markov process has been shown to approximate statistical properties of the original ergodic system $(\Omega, \mu, f)$ under a variety of different assumptions \cite{froyland2001extract}. We shall review some of them later in Section \ref{sec:connect:1}.   
	
	\paragraph{Goal} We have thus seen several aspects of the dynamical systems that go beyond its initial description as a map on an invariant space. The vast field of data-driven discovery of dynamical systems can be categorized based on which aspect is targeted. It is illuminating to consider the connections between these four aspects. In this section we have discussed how all the four aspects can be derived from Assumption \ref{A:f}. However, the connections are not clear from an approximation theory point of view. One may develop separate techniques to approximate $f$, the attractor $\Omega$, the Koopman operator $U$, and a Markov transition function $P$ separately. But that does not guarantee that the approximations are interchangeable. The precise conditions when they are interchangeable, and when they are not, are related to a deeper understanding of dynamical systems theory. Our goal is examine these inter-connections.
	
	\paragraph{Outline} While we shall have all the four aspects in mind, we begin with the primary one, i.e., the dynamics law $f$ itself. In Section \ref{sec:learn} we present a mathematical framework that unifies all the techniques for learning $f$. 
	We next consider the connections between the four aspects of dynamical systems, as present in Figure \ref{fig:outline1}. We survey the connections that are known to exist in Section \ref{sec:connect:1}. Of more interests are the lack of connections. We examine the challenges of approximating the invariant set from an approximation of $f$ in Section \ref{sec:stability}, and the challenges of approximating $U$ from an approximation of $f$ in  Section \ref{sec:predict}. The theoretical results in Section \ref{sec:predict} are illustrated by some numerical examples in Section \ref{sec:expt}. The discussion in the paper is summarized in Section \ref{sec:conclus}.
	
	\section{Learning the dynamics law} \label{sec:learn}
	
	We now fix our attention to the primary aspect of a dynamical system - its dynamics law $f$ as in Assumption \ref{A:f}. We describe an abstract framework presented in \cite{BerryDas_learning_2022} which unifies most of the techniques employed to learn $f$. It is formalized through the following Assumption which is built upon Assumptions \ref{A:f} and \ref{A:data} :
	
	\begin{Assumption} \label{A:pPhi}
		There is a map $\Phi : \Omega\to \real^L$ such that 
		\begin{enumerate} [(i)]
			\item $\Phi$ is an injective map;
			\item there is a function $g:\real^d\times \real^L \to \real^L$ such that $\Phi\circ f = g\circ \left(\phi \times \Phi\right)$.
		\end{enumerate}
	\end{Assumption} 
	
	The map $\phi$ is the measurement through which the dynamical system is observed. So the codomain of $\phi$ is often low dimensional, specially if $\phi$ is a partial observation. Since $\Phi$ is an injective map, it effectively serves as a representation of the dynamics-space $\Omega$ in the Euclidean space $\real^L$. Let $\dim_\mu$ denote the box-dimension of the invariant set supported by the measure $\mu$. We typically have
	\[ d \leq \dim_\mu << L . \]
	The task now is to reconstruct the dynamics using the maps $\phi, \Phi$. In a data-driven approach $\phi$ is assumed to be unknown. The function $\Phi$ could be explicit, such as in delay-coordinate techniques; or implicit, such as in invariant-graph based techniques.
	
	\paragraph{Embedding mechanism} The major abstraction provided by Assumption \ref{A:pPhi} is the function $g$. It represents the embedding mechanism underlying the learning technique, and is explicitly known and computable. Note that $\Phi\circ f$ is the evolution of $\Phi$ under one iteration of the dynamics of $f$. The relationship between $g, f$ and $\Phi$ is depicted in the following diagram :  
	\begin{equation} \label{eqn:paradigm:1}
		\begin{tikzcd}[row sep = large, column sep = large]
			\ &\real^d &\real^d\times \real^L \arrow[l, "g"'] \\
			\ &\Omega \arrow{u}{\Phi} &\Omega \arrow{u}{\phi\times\Phi} \arrow[l, "f"] \arrow{ul}{U\Phi}
		\end{tikzcd}
	\end{equation}
	The diagram in \eqref{eqn:paradigm:1} is called a \emph{commutative diagram}. It is a network in which each node is a space, and each directed edge a map. If two different paths connect the same pair of nodes, then the maps created by composing the maps along each path must be the same. These identities are known as \emph{commutation}-s. Commutation diagrams provide a holistic view of an arrangement of mathematical objects \cite{Shiebler2022kan, FongEtAl_backprop_2019, Das2023conditional, Das2024hmlgy}. At the same time it conveys in a pictorial but rigorous manner various identities that hold between the components. Here it depicts how $g$ stores the information of the map $f$ and embedding $\Phi$.
	
	\paragraph{Feedback function} Since $\Phi$ is an injective map (by Assumption \ref{A:pPhi}), the current state of $\Phi$ determines the current and all future states in $\Omega$ and therefore of $\phi$. Therefore for every $k\in\num$, there is a function $w_k$ such that
	\begin{equation} \label{eqn:def:wk}
		w_k:\real^L \to \real^d, \quad w_k \circ \Phi = U^k \phi = \phi \circ f^k. 
	\end{equation}
	The learning task can now be stated to be the determination of this map $w_k$, which we call the \emph{feedback function}. The most interesting case is when $k=1$, in which case we denote $w_k$ as $w$. Different learning techniques may have different $\Phi$ and $g$, but the goal is always to determine the $g$ which satisfies  \eqref{eqn:paradigm:1}.
	
	To better understand how the feedback function $w$ fits into our existing arrangement of Assumptions \ref{A:f}, \ref{A:data} and \ref{A:pPhi} we add some arrows to \eqref{eqn:paradigm:1} :
	\[\begin{tikzcd}[row sep = large, column sep = large]
		\ &\real^d &\real^d\times \real^L \arrow[bend right = 90]{dll}[swap]{\proj_1} \arrow{r}{\proj_2} \arrow{l}{g} & \real^L \\
		\real^d &\Omega \arrow{u}{\Phi} &\Omega \arrow{u}{\phi\times\Phi} \arrow{ur}{\Phi}   \arrow{ul}{U\Phi} \arrow{l}[swap]{f} \arrow[bend left = 30]{ll}[swap]{\phi}  
	\end{tikzcd}\]
	The arrows added are simply the projections from the Cartesian product space $\real^d\times \real^L$. The function $w$ can now be placed as a new branch in the diagram :
	\begin{equation} \label{eqn:paradigm:2}
		\begin{tikzcd}[row sep = large, column sep = large]
			\ &\real^d &\real^d\times \real^L \arrow[bend right = 90]{dll}{\proj_1} \arrow{r}{\proj_2} \arrow{l}{g} & \real^L \arrow{r}{w} &\real^d \\
			\real^d &\Omega \arrow{u}{\Phi} &\Omega \arrow{u}{\phi\times\Phi} \arrow{ur}{\Phi} \arrow{rr}[swap]{f} \arrow[bend right = 15]{urr}{U \phi} \arrow{ul}{U\Phi} \arrow{l}[swap]{f} \arrow[bend left = 30]{ll}[swap]{\phi} &\ &\Omega \arrow{u}[swap]{\phi} 
		\end{tikzcd}
	\end{equation}
	Our new additions have created two new commutations which turn out to be the identities 
	\[U\phi = \phi\circ f = w\circ \Phi . \]
	These identities follow from the very definitions of $U$ and $w$ respectively. Equation \eqref{eqn:paradigm:2} thus concisely summarizes the identities existing between the various spaces and functions discussed so far. This arrangement is present in all approaches to learning the dynamics-law $f$.
	
	\paragraph{The reconstructed system} Next note that the function $w_k$ in \eqref{eqn:def:wk} does not provide any information about the original manifold $\tilde\Omega$. The reconstruction of the dynamics and its law must be entirely in the Euclidean space $\real^d \times \real^L$. For $k=1$, such a dynamics is created by the following map :
	\begin{equation} \label{eqn:def:feedback_1}
		\reconstruct : \real^d \times \real^L \to \real^d \times \real^L, \quad 
		\left[ \begin{array}{c} u_{n+1} \\ y_{n+1} \end{array} \right]
		= \reconstruct \left[ \begin{array}{c} u_{n} \\ y_{n} \end{array} \right] 
		= \left[ \begin{array}{c} w \left( y_n \right) \\ g\left( u_{n}, y_n \right) \end{array} \right].
	\end{equation}
	The dynamics under the map $\reconstruct$ shadows the dynamics on $\Omega$. This is brought to light by redrawing \eqref{eqn:paradigm:2} :
	\begin{equation}\label{eqn:paradigm:3}
		\begin{tikzcd}[row sep = large, column sep = huge]
			\ & \blue{ \real^d\times\real^L } \arrow[blue, dashed]{r}{ \reconstruct } \arrow{dl}[swap]{\proj_2} & \blue{ \real^d\times\real^L}  \arrow{ddl}[pos=0.3]{\proj_1}\\
			\itranga{\real^L} \arrow[Itranga]{dr}[swap]{w} & \blue{\Omega}  \arrow[blue, dashed]{r}[pos=0.3]{f} \arrow[blue, dashed]{u}{\phi\times\Phi} \arrow[Itranga]{l}{\Phi} \arrow[Itranga]{d}[swap]{U\phi} \arrow{ur}{U\phi \times U\Phi} & \blue{\Omega}  \arrow{dl}{\phi} \arrow[blue, dashed]{u}[swap]{\phi\times\Phi} \\
			\ & \itranga{\real^d} &\ 
		\end{tikzcd}
	\end{equation}
	The commutation diagram in \eqref{eqn:paradigm:3} collectively presents the original dynamics $f$, the role of the embedding mechanism $g$ and the effect of the Koopman operator. The red loop presents the defining equation of the feedback function $w$. It conveys the same identity as \eqref{eqn:def:wk}. The top blue horizontal line represents the function
	\[ \reconstruct := \paran{ w\circ \proj_2 \,,\, g }  : \real^{d+L} \to \real^{d+L} \]
	The blue loop presents a \emph{semi-conjugacy} relation, 
	\[ \reconstruct \circ h = h \circ f, \quad h := (\phi, \Phi) : \Omega \to \real^{d+L} . \]
	Obtaining this relation is the entire purpose of the reconstruction procedure. Since the space $\tilde\Omega$ or the invariant set $\Omega$ are inaccessible, one can only to hope to recreate a dynamical system in data-space $\real^{d+L}$ that mimics the dynamics under $f$. The dynamics is via the map $\reconstruct$ and is tied to the dynamics of $f$ via the map $h$. The various colored loops in \eqref{eqn:paradigm:3} represent independent aspects of the learning problem and dynamics. They are tied together by the maps represented by black arrows. 
	
	The system \eqref{eqn:def:feedback_1} evolves in parallel to $f:\Omega\to\Omega$, while maintaining a conjugacy via the map $\phi\times\Phi$. The state vector $z_n := \left( u_n, y_n \right)$ is interpreted as follows : $u_n$ is represents the value $\phi(\omega_n)$, which is the unknown state $\omega_n$ observed through $\phi$. The point $y_n\in \real^L$ is the embedded point $\Phi(\omega_n)$ in $\real^L$. With this interpretation in mind, \eqref{eqn:def:feedback_1} is to be initialized as
	\begin{equation} \label{eqn:feedback_1_init}
		z_0 = z_0(\omega_0) := \left( \phi(\omega_0), \Phi(\omega_0) \right) \in \real^d \times \real^L .
	\end{equation}
	for some $\omega_0 \in \Omega$. This completes the description of the reconstruction. The perfect conjugation is possible only if the feedback $w$ is recovered perfectly. We next consider approximation errors in learning $w$.
	
	\paragraph{Hypothesis space} In a practical situation, $w$ or $w_k$ is estimated by looking for their closest approximation in a search-set or hypothesis space $\mathcal{H}$. The set $\mathcal{H}$ is typically a linear space, or a parameterized collection of functions. Thus the true function $w_k$ can be expressed as
	\begin{equation} \label{eqn:def:hat_w}
		w_k = \hat{w}_k + \Delta w_k ,
	\end{equation}
	where $\hat{w}_k \in \calH$ is the closest approximation, and $\Delta w_k$ is the error. Similarly to \eqref{eqn:def:feedback_1}, the dynamics under the approximated feedback function becomes
	\begin{equation} \label{eqn:def:feedback_2}
		\hat{\reconstruct} : \real^d \times \real^L \to \real^d \times \real^L, \quad 
		\left[ \begin{array}{c} u_{n+1} \\ y_{n+1} \end{array} \right]
		= \hat{\reconstruct} \left[ \begin{array}{c} u_{n} \\ y_{n} \end{array} \right] 
		= \left[ \begin{array}{c} \hat{w} \left( y_n \right) \\ g\left( u_{n}, y_n \right) \end{array} \right].
	\end{equation}
	We have already established that the task of discovering the dynamics law $f$ is essentially the task of approximating $w = w_1$. We shall later examine whether approximating $w$ is sufficient to approximate $U$. The difference between \eqref{eqn:def:feedback_1} and \eqref{eqn:def:feedback_2} will then bear significance. We show in  Section~\ref{sec:predict} that the asymptotic rate at which these two systems diverge could depend on both the spectral properties of the dynamics, or its Lyapunov exponents. 
	
	We do an extensive comparison of various learning techniques in Table \ref{tab:learn}. They employ different embedding mechanisms $g$ and different means of approximating $w$. There are mainly two paradigms for creating an embedding $\Phi$ - delay-coordinates and invariant-graph / echo-state networks (ESN) based techniques. The distinguishing features of these two techniques are summarized in Table~\ref{tab:prdgm}. We look at these two paradigms next.
	
	\begin{table}
		\caption{Various learning techniques for the feedback function, $w_k$ \eqref{eqn:def:wk}, all of which assume an embedding.}
		\begin{tabularx}{\linewidth}{L|L|L|L}
			Technique & Hypothesis space & Advantages & Disadvantages \\ 
			\hline
			Linear reg. & Linear combination of fixed basis functions or coordinates & Availability of techniques for linear cases & Poor fit for nonlinear functions \\
			\hline
			Kernel reg. \cite{BerryHarlim2017, BerryEtAl2015, DasDimitris_CascadeRKHS_2019} & $C^r(M)$ or $L^2(\mu)$ Spaces spanned by kernel sections or eigenvectors & Allows smooth interpolations and connections with underlying geometry  & Localized nature of basis functions require large number of basis functions \\
			\hline
			RKHS \cite{AlxndrGian2020, DasGiannakis_RKHS_2018, DasDimitris_CascadeRKHS_2019} & Span of eigenfunctions of kernel integral operators & Completely data-driven, allows out of sample extension & Inexplicit, unspecified nature of basis functions  \\
			\hline
			Nonlinear reg. \cite{Gallant1975} & Parameterized space of functions & Dependence on parameters allow application of manifold techniques such as gradient descent & Explicit knowledge of parameters as well as dependence on parameters required \\
			\hline
			Deep NNs \cite{LecunEtAl2015} & Functions parameterized by network activation and coupling parameters & Simplicity of implementation; scalability; explicit dependence on parameters known & Little apriori knowledge known about dimension of layers or number of layers; little knowledge about convergence rate of learning; huge number of variables to optimize \\
			\hline
			LSTM \cite{HarlimEtAl2021, MaEtAl_2018, Maulik_EtAl_2020} & Same as Deep NNs but with additional memory cells & Good for approximating functions which have sparse dependence over a long interval of time & Same as Deep NNs; more parameters to tune\\
			\hline
			Radial basis functions \cite{Smith1992} & Similar to kernel techniques & Provides a global representation of the map & Lack of normalization lead to non-uniformity in predictability\\
			\hline
			Local approximation techniques, such as simplex methods \cite{JimenezEtAl1992}, and local linear regression \cite{KugiumtzisEtAl98} & Nearest-neighbor based approximation of a neighborhood of the predictee point & Good approximation for low-curvature attractors, i.e., less oscillatory functions & Predictee point needs to be close to data cloud, feedback function unbounded.\\
			\hline
		\end{tabularx}
		\label{tab:learn}
	\end{table}
	
	\begin{table}
		\caption{The two learning paradigms satisfying Assumptions~\ref{A:f}, \ref{A:data} and \ref{A:pPhi} and the scheme in \eqref{eqn:paradigm:1}.}
		\begin{tabularx}{\columnwidth}{L|L|L|L}
			Name & $\Phi$ & Basis for convergence of $\Phi$ & $g$ \\ 
			\hline
			Invariant-graphs & Implicitly obtained \eqref{eqn:asdc9} & \eqref{eqn:cvd83} & Explicit : \eqref{eqn:cv93m} \\ 
			\hline
			Delay-coordinates & Explicitly obtained as basis functions & Ergodic convergence & Explicit : \eqref{eqn:sdjf3} \\
			\hline
		\end{tabularx}
		\label{tab:prdgm}
	\end{table} 
	
	\paragraph{Paradigm I : Invariant graphs} Techniques using this paradigm begin with a $C^1$-map $g:\real^d \times \real^L$ for which there is constant $\lambda\in (0,1)$ such that
	\begin{equation} \label{eqn:cv93m}
		\norm{ \nabla_y g(u, y) } \leq \lambda, \quad \forall u\in \real^d, \, y\in \real^L.
	\end{equation}
	A commonly used instance of $g$ is the function
	\[ g(x, y) = \tanh \paran{ Wy + b^Tx  } , \forall x\in \real^d, y\in \real^L , \]
	where $W$ is any matrix of norm $\leq 1$ and $b$ is any vector of length $d$. Using such a $g$ one can build a \emph{reservoir} system, which is a skew product system on $\Omega\times\real^L$ defined as
	\begin{equation} \label{eqn:def:reservoir}
		\left(\begin{array}{c} \omega_{n+1} \\ y_{n+1} \end{array}\right) := T_{\text{reservoir}} \left(\begin{array}{c} \omega_{n} \\ y_{n} \end{array}\right) := \left(\begin{array}{c} f(\omega_n) \\ g\left( \phi(\omega), y_n \right) \end{array}\right) .\end{equation}
	The paradigm of invariant graphs was first investigated in \cite{matthews1992uniform, jaeger2004harnessing}, and studied eventually in greater detail as \emph{echo-state networks} (ESN) \citep[e.g.][]{grigoryeva2014stoch, GrigoryevaHartOrtega2021chaos, LukoJaeger2009,  GononOrtegafading_fading_2021, GrigoryevaHartOrtega2021}, and \emph{reservoir computers} \cite{grigoryeva2021learn, LuEtAL2017}.  Although \eqref{eqn:def:reservoir} involves the underlying dynamical map $f$, the actual knowledge of $f$ is not needed. Note that the dynamics in the $y$-coordinate is linked to the $\omega$-coordinate through the measurement $\phi$. 
	
	ESNs are used in two phases. In the training phase, one provides as input the measurements $\left\{ \phi(\omega_n) \right\}_{n=0}^{N}$. Thus $(\Omega, f)$ remains unknown but continues to drive the reservoir variable $y$. The variable $y$ settles down into a representation of the attractor, in a manner which we make precise below.
	
	\begin{proposition} \label{thm:kncd3l}
		Let Assumption \ref{A:f} hold, and $\phi:\Omega\to \real^d$ be a continuous map, and $g:\real^d\times \real^L \to \real^L$ be a $C^1$ map satisfying \eqref{eqn:cv93m}. Then
		\begin{enumerate}[(i)]
			\item there is a map $\Phi:\Omega \to \real^L$ such that Assumption \ref{A:pPhi}~(ii) is satisfied.
			\item The graph of $\Phi$ is invariant, i.e.,
			\begin{equation} \label{eqn:asdc9}
				\left( U^n\Phi \right)(\omega) := \Phi \left( f^n(\omega) \right) = \proj_Y T_{\text{reservoir}}^n \left( \omega, \Phi(\omega) \right) , \quad \forall n\in\num , \quad \forall \omega\in \Omega.
			\end{equation}
			\item The graph of $\Phi$ is globally attracting, i.e.,
			\begin{equation} \label{eqn:cvd83}
				\lim_{n\to\infty} \left( U^n\Phi \right)(\omega) = \lim_{n\to\infty} \proj_Y T^n\left( x, y \right), \quad \forall x\in \bar{\Omega}, \, y\in Y.
			\end{equation}
			\item If $\norm{ \frac{\partial }{\partial_u} g } \leq 1$, then  Assumption~\ref{A:g_contract} is also satisfied.
		\end{enumerate}
	\end{proposition}
	
	Assumption~\ref{A:g_contract} is an additional assumption requiring that $g$ be a non-expansive map in each of its variables. It is described later Section~\ref{sec:stability} and is used to establish stability properties.  Parts (i)---(iii) are immediate consequences of results by Stark \cite{Stark1999}, or by Grigoryeva et. al. \citep[][Thm III.1]{GrigoryevaHartOrtega2021}. 
	
	The invariant graph property leads to a fulfillment of the identity in Assumption \ref{A:pPhi}~(ii). However, the injectivity condition of Assumption \ref{A:pPhi}~(i) remains to be proven rigorously. It has been generally observed that for $L$ large enough, $\Phi$ is also injective. The ground Assumption~\ref{A:f} is assumed while running the system. Note that the map $\Phi$ is not obtained explicitly but implicitly through the state variables of the network. Given any arbitrary initialization to \eqref{eqn:def:reservoir}, by \eqref{eqn:cvd83}, the internal states of the reservoir converge to an invariant graph over $\Omega$. The function $\Phi$ is precisely the function whose graph is invariant. Although it will remain indeterminate, its values over a dynamic trajectory, i.e., the values $\phi(f^n \omega_0)$ will be obtained, for some unknown initial point $\omega_0$. Thus the family of ESNs fall within the framework of Assumptions \ref{A:f}, \ref{A:data} and \ref{A:pPhi}. 
	It is still an open question whether the map $\Phi$ representing the invariant graph is injective. Some results are known for special cases \cite{GrigoryevaHartOrtega2021chaos, HaraKokubu2024learn}. 
	
	\paragraph{II. Delay coordinates} An effective and numerically inexpensive means of obtaining an embedding of a dynamical system is using delay coordinates \cite{DasGiannakis_delay_2019, BerryHarlim2017, SauerEtAl1991, BerryEtAL2013}. To relate to our framework, fix a number of delays $Q\in\num$, and set 
	%
	\begin{align} \label{eqn:sdjf3} 
		L = Qd, &\; \\
		\Phi : \Omega\to \real^L, &\; \omega \mapsto \paran{ \phi\left( \omega \right), \ldots, \phi\left( f^{Q-1} \omega \right) } \\
		g : \real^d\times \real^L \to \real^L, &\; \paran{ u, y^{(1)}, \ldots, y^{(Q)} } \mapsto \paran{ u, y^{(1)}, \ldots, y^{(Q-1)}  }
	\end{align}
	%
	Thus Assumptions \ref{A:f}, \ref{A:data} and \ref{A:pPhi} are already satisfied. Now \eqref{eqn:def:feedback_1} becomes :
	\[ \reconstruct_{\text{delay-coord}} :\real^{d\times dQ} \to \real^{d\times dQ} := \left[\begin{array}{c}
		u \\
		y^{(1)} \\
		y^{(2)} \\
		\vdots \\
		y^{(Q)} 
	\end{array}\right] \mapsto  
	\left[\begin{array}{c}
		w\left( y^{(1)}, \ldots, y^{(Q)} \right) \\
		y^{(1)} \\
		\vdots \\
		y^{(Q-1)} 
	\end{array}\right]
	\]
	Note that in this case, $g$ is a linear map. We have :
	
	\begin{proposition} \label{prop:dwa4ep}
		\cite{SauerEtAl1991} Let Assumption~\ref{A:f} hold. Then for a typical map $\phi:\Omega\to \real^d$, if $Q\in\num$ is large enough, then $\Phi$ defined through \eqref{eqn:sdjf3} is an injective, and thus Assumption~\ref{A:pPhi} is satisfied. Moreover, Assumption~\ref{A:g_contract} is also satisfied.
	\end{proposition}
	
	This general result establishes the injectivity condition of Assumption \ref{A:pPhi}~(i) in the typical case. Proposition \ref{prop:dwa4ep} provides the mathematical guarantee on the reliability of the delay-coordinate embedding method. A vast number of timeseries analysis techniques have implicitly or explicitly relied on this result.   
	
	This completes our demonstration of how the framework depicted in \eqref{eqn:paradigm:3} unifies all the common techniques of learning the dynamics law. Some interesting connections from the view point  of \emph{generalized synchronization} may be found in \cite{HartHookDawes2020embed, HartHookDawest2021Tikh}.
	
	\section{Connecting different descriptions} \label{sec:connect:1} 
	
	\begin{figure}
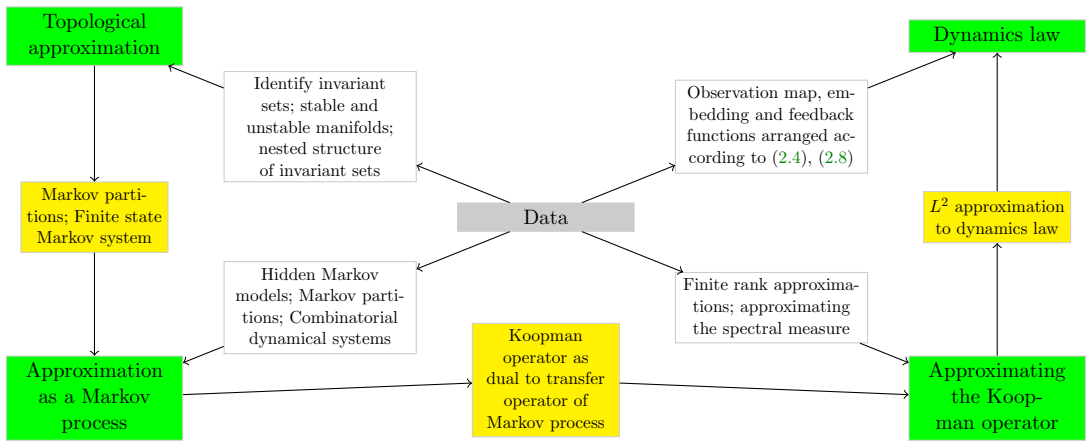
\center
		\begin{tikzpicture}[scale=0.6, transform shape]
			\input{outline_L1.tex}
			\input{outline_L2.tex}
		\end{tikzpicture}
		\caption{Relating different aspects of a dynamical system. Figure \ref{fig:outline1} has been expanded by adding connections (yellow) between the four aspects (green) of dynamical systems. }
		\label{fig:outline2}
	\end{figure} 
	
	\begin{figure}
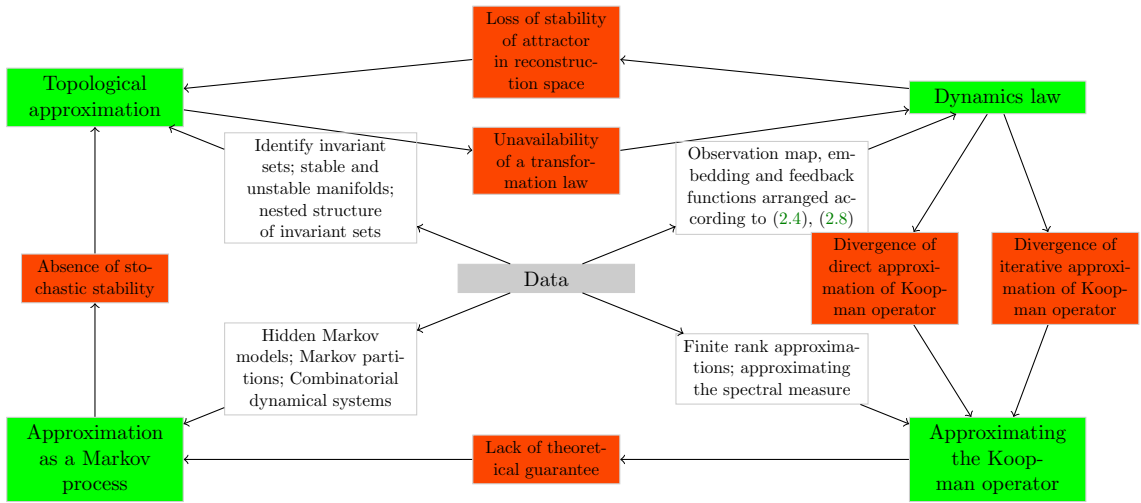
\center
		\begin{tikzpicture}[scale=0.6, transform shape]
			\input{outline_L1.tex}
			\input{outline_L3.tex}
		\end{tikzpicture}
		\caption{Gaps in learning techniques for dynamics. Figure \ref{fig:outline1} has been expanded by adding missing links (red) that should connect the four aspects (green) of dynamical systems.}
		\label{fig:outline3}
	\end{figure} 
	
	\begin{figure}\center
		\includegraphics[width=.95\linewidth]{\figs 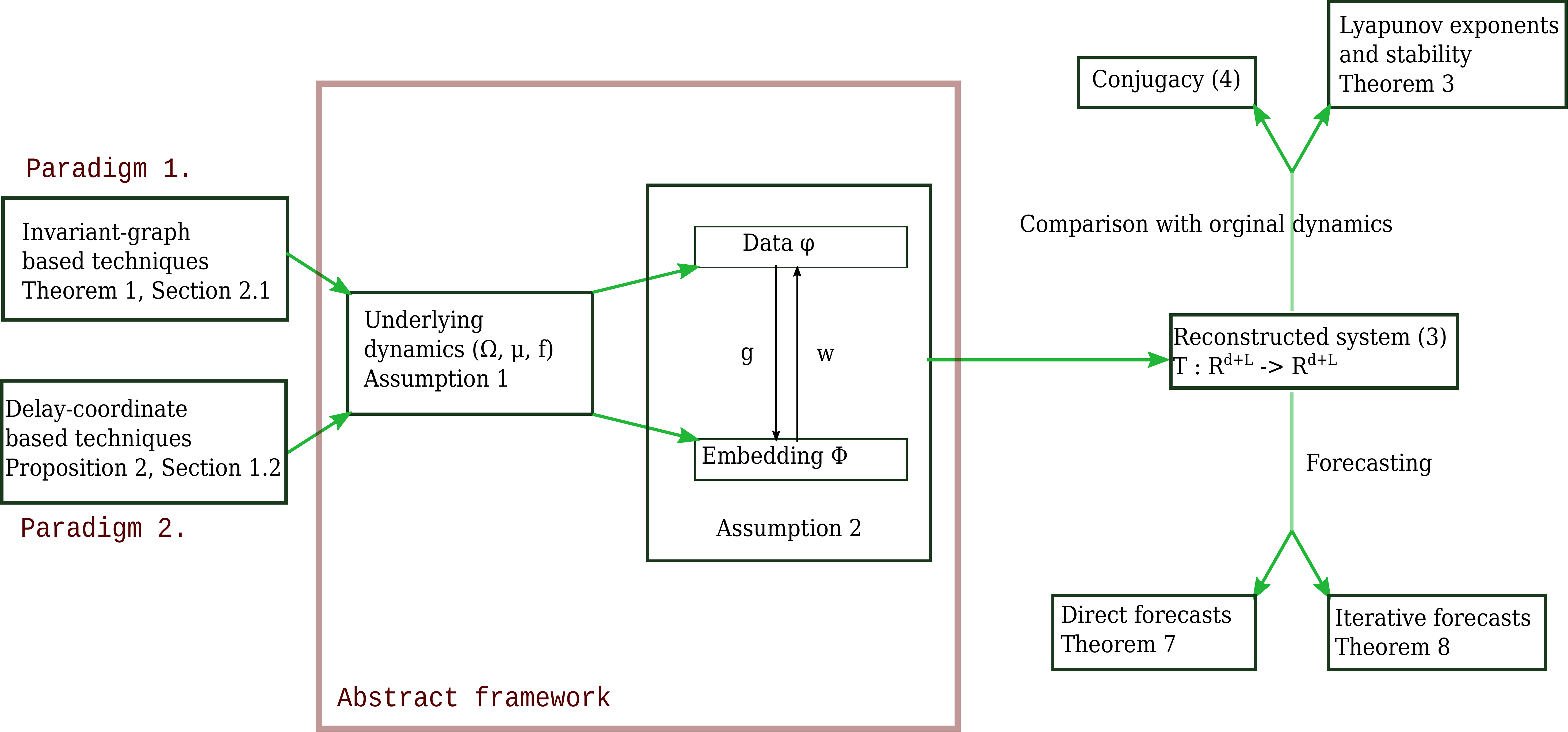}
		\caption{Learning the dynamics law. Assumption~\ref{A:pPhi} provides a unifying framework for different techniques of learning the dynamics law $f$. The details of this framework is explained in detail in Section \ref{sec:learn} and \eqref{eqn:paradigm:1}. The two main paradigms of learning dynamics - based on invariant graphs and delay coordinates respectively, are particular instances of this framework. Sections \ref{sec:stability} and \ref{sec:predict} examine whether approximating the dynamics law is sufficient to approximate the Koopman operator and topological attractor. These aspects of the dynamics have been depicted as the issues of predictability and stability of the reconstructed system respectively. and in which the embedding is implicit and explicit respectively. }
		\label{fig:outline4}
	\end{figure}
	
	We have seen the four aspects of dynamical systems in Sections \ref{sec:intro} and reviewed various means by which they are reconstructed from data. The primary aspect is the dynamics law from which each of the other three can be derived. We have focused on this primary aspect in Section \ref{sec:learn}. At this stage we consider the question of whether approximating any one aspect leads to an approximation of another aspect. Figure \ref{fig:outline2} presents some such connections. These connections are based on deeper theoretical understanding of dynamical systems, and we shall take a quick look at them.
	
	\paragraph{ Connection 1} We have already discussed how a topological approximation of the invariant set could lead to an approximation of the dynamics as a Markov process. A topological approximation of the invariant set leads is effectively an approximation of its neighborhood. Such a neighborhood can be partitioned into $m$ cells, followed by a computation of the transition matrix \eqref{eqn:def:transition}. The matrix simulates a Markov process on the finite state space $\braces{1, \ldots, m}$ which could serve as a discretized Markov approximation of the system. 
	
	As mentioned before, there are a variety of conditions on the dynamics $(f,\Omega,\mu)$ under which this finite state Markov process mimics the ergodic behavior of $\mu$. A condition typically relied on is quasi-compactness \cite[see]{HuntMiller1992approx} of the transfer operator $\calP$. The latter is a property borne by the spectrum of $\calP$. If $\calP$ is quasi-compact, then it has a unique fixed point. Moreover all ergodic sums converge to the unique fixed point \cite[Thm 2.3]{DingLiZhou2002fin}. Techniques such as the simplicial triangulation method \cite{DingLiZhou2002fin} are shown to provide good convergence to the fixed stationary density. Another important property concerns the iterations of $\calP$ and is said to be of \emph{Lasota Yorke type} \cite[e.g.]{DingDuLi1993high}. Simple techniques such as continuous piecewise linear approximations are shown to be convergent \cite[see]{DingZhou2001constr}.
	Some other conditions are unique ergodicity \cite[e.g.]{Hunt1998unique} a condition on the global stability of the $\calP$; and algebraic conditions on the transition matrices \citep[e.g.]{Froyland1998approx}.
	
	All these conditions are about the ergodic system $(\Omega, f, \mu)$. Thus under any of these conditions, a topological approximation does guarantee an approximation as a Markov process.
	
	\paragraph{ Connection 2} We now examine how an approximation of the dynamics as a Markov process could lead to an approximation of the Koopman operator $U$. Given any stationary stochastic process $X_t$ one can define the transition probability as a function which which assigns to every time $t$ and point $x$ in the state-space the probability measure 
	\[ P(t, x)(A) := \mathbb{E} \paran{ X_t\in A \,|\, X_0 = x } . \]
	The Koopman operator $U^t$ \cite{Philipp2024error} on $L^\infty(\mu)$ can be simply defined as
	\[ \paran{ U^t \psi }(x) := \mathbb{E} \paran{ \psi(X_t) \,|\, X_0 = x } . \]
	Alternatively one can also define it directly in terms of the transition probability as
	\[ \paran{ U^t \psi }(x) = \int \psi dP(t,x) . \]
	This is the stochastic analog of the definition \eqref{eqn:def:Koop}. A more natural interpretation of the Koopman operator for stochastic processes is as a dual operator to the transfer operator \cite{cvitanovic1998trace, schutte2001transfer}. Any measure $\nu$ on the state space $\calX$ transforms under the stochastic process as
	\[ P^t \nu := \int P(t,x) d\nu(x) . \]
	The evolution of the measure $\nu$ is linear w.r.t. $\nu$ and is tracked through the transfer operator $P^t$. A stationary measure for the process is a measure $\mu$ such that $P^t \mu = \mu$. The transfer operator is one of the many ways of tracking the propagation of probability measures on the state space. Now suppose that a reference probability measure $\nu_0$ is fixed, and there is a probability measure $\nu$ which is absolutely continuous w.r.t. $\nu_0$. So $\nu$ has a density function $\rho \in L^1(\nu_0)$. Then for all times $t$, the measure $P^t \nu$ is also absolutely continuous w.r.t. $\nu_0$. So it has a density function $\rho_t \in L^1(\nu_0)$. The dependence of $\rho_t$ on $\rho$ is also linear, leading to an interpretation of $P^t$ as an operator on $L^1(\nu_0)$. Now suppose $\nu_0=\mu$, a stationary measure. Then the Koopman operator turns out to be the dual to the transfer operator :
	\[ \left\langle P^t \alpha, \beta \right\rangle_{L^2(\mu)}  = \left\langle \alpha, U^t \beta \right\rangle_{L^2(\mu)} , \quad \forall \alpha\in L^1(\mu) , \, \beta \in L^\infty(\mu) . \]
	This weak characterization of the Koopman operators have been utilized in many different ways \cite{KlusSchutte2016, KlusEtAl2017, Philipp2024error} for a numerical approximation of the Koopman operator.
	
	\paragraph{ Connection 3} We now examine how an approximation of the Koopman operator $U$ could lead to an approximation of the dynamics law $f$ itself. Since $U$ is an infinite dimensional operator, its approximation usually involves the choice of a finite dimensional subspace $\calV$. For simplicity we assume that the approximation is of the form
	\[ \hat{U} := \proj_{\calW} U \proj_{\calW} . \]
	Now let $\gamma = \Id_{\real^D}$, the identity function on the reconstruction space $\real^D$. Then note that
	\[ U \gamma = \gamma \circ F = F . \]
	Thus the reconstructed dynamics law $F$ is the image of the identity function under the Koopman operator. The law can be approximated as
	\[ F \approx \hat{F} := \hat{U} \gamma =  \proj_{\calW} U \proj_{\calW} \paran{ \Id_{\real^D} } . \]
	As the size of the approximation space $\calW$ increases, $\hat{F}$ converges in some appropriately chosen norm, to $F$. The approximation of the dynamics law from the Koopman operator follows from the fact that the Koopman dynamics contains the original phase space dynamics within its description.
	
	This completes a discussion of all the connections shown in yellow in Figure \ref{fig:outline2}. These connections are directional and not mutual. In fact many of the connections cannot be reversed using prevailing techniques in dynamical systems theory. These are outlined in red in Figure \ref{fig:outline3}. We review some of these next.
	
	\paragraph{ Dis-connection 1} An approximation as a Markov process does not necessarily provide an approximation of the invariant set. In a Markov approximation, the deterministic iterations of the dynamics is represented by Markov transitions, which are not deterministic but have a distribution. From a numerical point of view, the aspect of a Markov process visible to experiments is its stationary measure $\nu$. Thus in a good Markov approximation one can expect $\nu$ to be a good statistical approximation of $\mu$. However, one cannot in general guarantee that diminishing the spread or entropy of the Markov transition function would guarantee the convergence of $\nu$ to $\mu$. Ergodic systems which have this property are called \emph{stochastically stable}. If the ergodic system $(\Omega, \mu, f)$ has a sequence of Markov transitions in which the distributions of the Markov transitions converge to a Dirac delta function, then it is called a \emph{zero-noise} limit. Stochastic stability means that zero-noise limits also guarantee an approximation of the measure. Stochastic stability is difficult to ascertain, and only been proved for special systems \cite{kushner2006stoch, Young_StochHyp_1986, CowiesonYoung2005, meyn2012markov}. In fact even in low dimensional deterministic systems, invariant regions are shown to have extremely sensitive and even dis-continuous dependence on the map \cite[e.g.]{DasYorke2020, DasSaddles2015, DasJim17_chaos}. A recent technique \cite{Das2024zero} creates a sequence of Markov processes such that the dynamics along with its invariant set is the zero-noise limit of these processes. Thus the invariant measure $\mu$ may not be approximable, but its support $\Omega$ may be.
	
	\paragraph{ Dis-connection 2} An approximation of the Koopman operator $U$ can be utilized to derive a Markov process, but it remains hard to decide whether this process approximates the actual dynamics. Suppose that $\calV = \braces{V_1, \ldots, m}$ is a partitioning of some small neighborhood of $\Omega$. A Markov approximation involves computing the $m\times m$ transition matrix $\mathbb{P}_{i,j} := \mu \paran{ V_i \cap f^{-1} V_j }$. Now note that 
	\begin{align*}
		\mathbb{P}_{i,j} &:= \mu \paran{ V_i \cap f^{-1} V_j } = \int_{V_i} \paran{ U 1_{V_j} } = \int_{\Omega} 1_{V_i} \paran{ U 1_{V_j} } d\mu
		&= \bracketBig{ 1_{V_i} , U V_j }_{L^2(\mu)} = \bracketBig{ \calP 1_{V_i} , V_j }_{L^2(\mu)} .
	\end{align*}
	where $\calP$ is the transfer operator. Thus from any approximation $\hat{U}$ of $U$ we get the approximated matrix :
	\[ \hat{ \mathbb{P} }_{i,j} = \bracketBig{ 1_{V_i} , U V_j }_{L^2(\mu)} \approx \bracketBig{ 1_{V_i} , \hat{U} 1_{V_j} }_{L^2(\mu)} . \]
	This matrix $\hat{ \mathbb{P} }$ can be column normalized so that it mimics the generator of Markov process of the finite state space $\braces{1, \ldots, m}$. 
	This simple principle connects two completely different entities. The Markov process is stochastic, nonlinear and on a finite state space. On the other hand $U$ is linear, deterministic, and acts on an infinite dimensional functional space. 
	
	In spite of the simplicity of the process, we must recall that an arbitrary partition based matrix $\mathbb{P}$ may not provide an approximation of the dynamics under $f$ as a Markov process. Although we have discussed some conditions for $(\Omega, \mu, f)$ there are not many verifiable conditions on $U$ which guarantee a Markov approximation. Thus lack of theoretical guarantee this counts as a disconnection.
	
	Figure \ref{fig:outline3} presents an overview of several such lack of connections. We next look in more depth at the disconnects from the dynamics law to the other aspects of the dynamics. These are laid out more lucidly in Figure \ref{fig:outline4}. The issues presented on the right side of this figure are investigated in turn in Sections \ref{sec:stability} and \ref{sec:predict} next. 
	
	\section{Stability of reconstructed system} \label{sec:stability} 
	
	We now examine whether an approximation of the invariant set is possible from an approximation of the dynamics law. At this stage, the issue of visibility discussed in Section \ref{sec:intro} becomes important. The learning framework presented in Section \ref{sec:learn} indicates that the original manifold $\tilde\Omega$ is out of reach, and the learning and reconstruction must take place in different Euclidean spaces. The image of $ h := \phi\times \Phi$ is a bijective image of $\Omega$, and is invariant under the dynamics of $\reconstruct$ \eqref{eqn:def:feedback_1}. Thus irrespective of whether $\Omega$ or $\mu$ is visible in $\tilde\Omega$ to $f$, $h(\Omega)$ or $h_*\mu$ may not be visible to $\reconstruct$ in $\real^{d+L}$. The question becomes even more difficult if one also takes into consideration approximation errors in estimating $F$. 
	
	The issue of visibility can be interpreted as the issue of stability and/or instability. It is quantified best through global indexes called Lyapunov exponents. One needs to quantify the rate of deviation of perturbations from the set $X := h(\Omega)$ under the map $\reconstruct$ Let the distinct Lyapunov exponents of $(\Omega, f, \mu)$ be $\lambda_1 > \lambda_2 \cdots > \lambda_r$, with corresponding Oseledets splitting $T\Omega = E_1 \oplus \cdots \oplus E_r$. Since the dimension of $\tilde\Omega$ is $m$, the multiplicities of the $\lambda_i$ sum to $m$. The $E_i$s corresponding to negative valued $\lambda_i$ constitute the stable directions for $\Omega$ in $\tilde\Omega$. Similarly the $E_i$ corresponding to positive valued $\lambda_i$ constitute the unstable directions. Moreover,
	\[ \lim_{n\to\infty} \frac{1}{n} \ln \norm{Df^n(\omega) v_i} = \lambda_i, \quad \mu-a.e. \omega\in \Omega, \, \forall v_i\in E_i(\omega)\setminus\{0\} . \]
	The map $\reconstruct$ acts in the higher dimensional ambient space $\real^{L+d}$ and therefore will have $L+d$ Lyapunov exponents (counting multiplicities). We show in Theorem~\ref{thm:lambda1}(i) that $m$ of the $d+L$ Lyapunov exponents of $\reconstruct$ coincide with the original $m$ Lyapunov exponents of $f$. We are interested in these other $d+L-m$ Lyapunov exponents of the reconstructed systems, and their positions relative to $\lambda_1(f), \ldots, \lambda_d(f)$. It is an extremely challenging task to estimate or even control the magnitude of the additional Lyapunov exponents. 
	
	By the very definition of Lyapunov exponents, the $\lambda_i(\reconstruct)$ depend not only on the invariant set $X$ but also on its neighborhood.An essential part of $\reconstruct$ is the feedback function $w$. The function $w : \real^L \to \real^d$ is defined uniquely only on $X$. The conjugacy in \eqref{eqn:paradigm:3} will be preserved on $\Omega$ irrespective of the nature of the extension of $w$ to a neighborhood of $X$. We define a collection 
	\[ \W := \SetDef{ \hat{w} \in C^1\left( \real^L; \real^d \right) }{  \hat{w} \rvert_{X} = w\rvert_{X} }, \]
	equipped with the $C^1$-topology. Every $\hat{w}\in \W$ is a $C^1$ function satisfying $\hat{w}\circ\Phi(\omega) = (U\phi)(\omega)$ for every $\omega\in \Omega$. Any choice of $\hat{w}\in \W$ leads to a different dynamics in $\real^{L+d}$ as in \eqref{eqn:def:feedback_1}. However $X$ continues to remain an invariant ergodic set. Therefore the top Lyapunov exponent $\lambda_1$ of $\reconstruct|X$ will be a function of $\hat{w}$. This leads to the function :
	\begin{equation} \label{eqn:def:lambda1}
		\lambda_1 : \W\to \real, \quad \lambda_1(\bar{w}) := \lambda_1( \reconstruct ).
	\end{equation}
	An important realization in the learning theory for dynamics \cite{BerryDas_learning_2022} is that the stability of the reconstruction is directly captured via \eqref{eqn:def:lambda1}. Our goal will be to study how close $\lambda_1(\reconstruct)$ can be made to $\lambda_1(f,\mu)$.
	
	\paragraph{Stability gap} As pointed out in \cite{DechertGencay2000, DechertGencay1996}, the top Lyapunov exponent of the reconstructed system may exceed that of the original system. Moreover, some of the additional Lyapunov exponents may be positive. All of these contribute to additional instabilities being introduced into the system. We define the \emph{stability gap} of the reconstruction as
	\[ \text{stability gap} :=\inf_{\bar{w}\in \W} \lambda_1(\bar{w}) - \lambda_1(f,\mu) . \]
	The stability gap is always non-negative, as will be shown in Theorem~\ref{thm:lambda1}. We shall study ways to obtain a bound on the stability gap in our next theorem. We shall need two additional assumptions. This first is about the embedding mechanism :
	
	\begin{Assumption} \label{A:g_contract}
		The function $g$ from Assumption~\ref{A:pPhi} further satisfies :
		\[ \sup_{\omega\in\Omega} \norm{\partial_1 g} \rvert_{(\phi(\omega), \Phi(\omega))} \leq 1 , \quad \sup_{\omega\in\Omega} \norm{\partial_2 g} \rvert_{(\phi(\omega), \Phi(\omega))} \leq 1 . \]
	\end{Assumption}
	
	To get more control of the neighborhood of $X$ we shall bring in an assumption about having a retraction to the range of $\Phi$. Let $\mathcal{U}$ be a neighborhood of $\ran \Phi$. Recall that a continuous map $\retract : \mathcal{U} \to \ran \Phi$ is said to be a \emph{retract} if $\retract|_{ \ran \Phi } = \Id_{ \ran \Phi }$.
	
	\begin{Assumption} \label{A:w_ext}
		There is a continuous retraction $\retract : \mathcal{U} \to \ran \Phi$, for some open neighborhood $\mathcal{U}$ of $\ran \Phi$ in $\real^L$.
	\end{Assumption}
	
	A retraction of a set is a continuous deformation of its neighborhood into itself. All submanifolds have neighborhood retraction within their ambient manifold. However, it is in general difficult to ascertain when a retraction exists. For the retraction $\retract$, we are interested in the Lipschitz norm of the retraction
	\begin{equation} \label{eqn:def:kappa_ret}
		\kappa_{\retract} := \sup_{y\in \ran\Phi} \limsup_{y'\to y} \frac{ d\left( \retract(y), \retract(y') \right) }{ d\left( y, y'\right) } .
	\end{equation}
	Finally, we identify a quantity $\Csens$ as the \emph{sensitivity constant} that depends on the (fixed) functions $\phi, \Phi$ and a point $\omega\in \Omega$.
	\[ \Csens : \Omega \to \real^+, \quad \Csens(\omega) := \sup \SetDef{ \frac{ \norm{D\phi(\omega) v } }{ \norm{D\Phi(\omega) v } } }{ v\in T_\omega \Omega \setminus \{0\} }. \]
	The quantity $\Csens$ measures how sensitive the measurement map $\phi$ is towards changes in the embedding coordinates $\Phi$. Both $\retract$ and $\Csens$ will be used to bound the gap between between $\lambda_1(\reconstruct)$ and $\lambda(f)$.
	
	\begin{theorem}[Stability of reconstruction] \label{thm:lambda1}
		\cite[Thm 3.1]{BerryDas_learning_2022}
		Let Assumptions~\ref{A:f}, \ref{A:data} and \ref{A:pPhi} hold. Then
		\begin{enumerate} [(i)]
			\item The $d+L$ Lyapunov exponents of $\reconstruct$ contains as a subset the $m$ Lyapunov exponents of $f$.
			\item $\lambda_1(\bar{w})$ is upper semi-continuous with respect to $\bar{w}$. In other words, for every $\epsilon>0$, there is a $C^1$ neighborhood $\mathcal{U}$ of $\bar{w}$ such that
			\[ \lambda_1( \bar{w'} ) < \lambda_1(\bar{w}) + \epsilon, \quad \forall \bar{w'} \in \mathcal{U} . \]
			\item Suppose Assumptions~\ref{A:g_contract} and \ref{A:w_ext} also hold. Then the stability gap is bounded by
			\begin{equation} \label{eqn:smf03}
				\inf_{\bar{w}\in \W} \lambda_1(\bar{w}) - \lambda_1(f,\mu) \leq \int \ln \left[ 1 + \left(1 + \Csens(\omega) \right) \kappa_{\retract} \right] d\mu(\omega) .
			\end{equation}
		\end{enumerate}
	\end{theorem} 
	
	Claims (i) and (ii) of Theorem~\ref{thm:lambda1} are immediate consequences of results from \cite{DechertGencay1996, BochiViana_conti_2005, Viana_Lyap_2020}. Assumption~\ref{A:w_ext} is of a topological nature and depends on the topological or geometrical properties of $X$. Equation \eqref{eqn:smf03} gives a general global bound on the stability gap. In practice, $\hat{w}$ is obtained from some hypothesis space which is determined by the application domain. In such situations there is no guarantee of the stability being preserved up to an $\epsilon$ error. The following corollary applies to the use of a large number of delay-coordinates.
	
	\begin{corollary} \label{cor:delay_stab}
		Let $\Psi^t : \Omega\to \Omega$ be a smooth flow and $f$ be the time-$\Delta t$ map $f=\Psi^{\Delta t}$. Let all the conditions in Assumptions \ref{A:f}, \ref{A:data} and \ref{A:pPhi} be met and the delay coordinate paradigm \eqref{eqn:sdjf3} be implemented. Suppose further that there is a retraction map as in Assumption~\ref{A:w_ext} for which the Lipschitz constant $\kappa_{\retract} = 1$. Then there is a constant $C_2$ depending only on the flow such that $\Csens(\omega) \leq \frac{1}{Q} + 0.5 C_2 Q \Delta t$ for every $\omega\in \Omega$. In particular, 
		\[ 0\leq \inf_{\bar{w}\in \W} \lambda_1(\bar{w}) - \lambda_1(f,\mu) \leq \ln\left[ 2+ \frac{2}{Q} + C_2 Q\Delta t \right]. \]
	\end{corollary}
	
	The criterion that $\kappa_{\retract} = 1$ is attained for example when $X = \ran\Phi$ is a manifold, and $\retract$ is a tubular neighborhood retract. 
	
	Even if the reconstruction has a zero stability gap, it does not guarantee that $X$ would be visible. A nonzero stability gap indicates the presence of new directions of instability. At present there is no known technique of learning the dynamics law that would simultaneously achieve a zero or small stability gap. This is precisely one of the challenges presented in Figures \ref{fig:outline3} and \ref{fig:outline4}.
	
	\section{Forecasts with reconstructed system} \label{sec:predict} 
	
	We now turn our attention to the question of whether an approximation of the dynamics law $f$ is enough to guarantee an approximation of the Koopman operator $U$. It is convenient to interpret $U$ as a linear representation of forecasting. Thus the discussion in this section may be interpreted as an evaluation of the effectiveness of reconstruction models as tools of forecast / prediction. From the very outset, we make a distinction between two types of forecasting, \emph{direct} or \emph{iterative}, as suggested by Casdagli \cite{Casdagli1989}. They are respectively the tasks of estimating the $k$-th iterate $U^k$ directly, and iterating $k$-times an estimate of the base operator $U$. 
	
	We now state this more precisely in the context of Assumption \ref{A:f} and \ref{A:data}. Fix an $\omega_0 \in \Omega$, and let $z_0 = h(\omega_0) \in \real^{d+L}$. There are two ways of estimating the value $\phi(f^k \omega_0)$ after $k$ iterations of the base dynamics : we can iterate \eqref{eqn:def:feedback_2} $k$ times, and the first coordinate of $\hat{\reconstruct}^k z_0$ will serve as an approximation of $\phi( f^k \omega_0) = (U^k \phi)(\omega_0)$. We call this the iterative method, and its accuracy can be estimated via
	\begin{equation} \label{eqn:def:err_iter}
		\begin{split}
			\text{error}_{\text{iter}}(k, \omega) &:= \norm{ U^k \phi (\omega) - \proj_1 \circ \hat{\reconstruct}^k \circ (\phi, \Phi)(\omega) }_{\real^d} , \\
			\text{error}_{\text{iter}}(k) &:= \left[ \int_\Omega \text{error}_{\text{iter}}(k, \omega)^2 d\mu(\omega) \right]^{1/2}.
		\end{split}
	\end{equation}
	Alternatively one can directly approximate $w_k$ via \eqref{eqn:def:hat_w} and obtain a \emph{direct} estimate. The corresponding errors are
	\begin{equation} \label{eqn:def:err_direct}
		\begin{split}
			\text{error}_{\text{direct}}(k, \omega) &:= \norm{ U^k \phi (\omega) - \hat{w}_k \circ \Phi (\omega) }_{\real^L} ,\\
			\text{error}_{\text{direct}}(k) &:= \norm{ U^k \phi - \hat{w}_k \circ \Phi }_{L^2(\mu)} = \left[ \int \text{error}_{\text{direct}}^2 (k, \omega) d\mu(\omega) \right]^{1/2}.
		\end{split}
	\end{equation}
	To aid the discussion, we briefly review some spectral aspects of the Koopman operator.
	
	\paragraph{Spectral aspects} The operator $U$ and all its iterates $U^n$ are unitary maps. As a result their spectrum lies on the unit circle in the complex plane. The constant function $1_\Omega$ is always an eigenfunction for $U$, with eigenvalue $1$. A consequence of $\mu$ being ergodic is that $1$ is a simple eigenvalue. Let $\Disc$ denote the closure in $L^2(\mu)$ of all eigenfunctions of $U$. Then one has the orthogonal decomposition
	\begin{equation} \label{eqn:def:L2split}
		L^2(\mu) = \Disc \oplus \Disc^\bot .
	\end{equation}
	These two components $\Disc, \Disc^\bot$ not only have different ergodic properties \cite{Halmos1956, DasJim2017_SuperC, DSSY2017_QQ, DSSY_Mes_QuasiP_2016}, but also respond differently to data-analytic and harmonic analytic tools \cite{DasGiannakis_delay_2019, DasGiannakis_RKHS_2018}. Eigenfunctions reveal non-mixing, non-chaotic dynamics embedded within $(\Omega,f)$. The related phenomenon of embedded Lie-group rotations \cite{Das2023Lie} is also of similar interest. 
	
	One of the major consequences of the splitting \eqref{eqn:def:L2split} is seen in the decay of correlations :
	\begin{equation} \label{eqn:untry:2}
		\bracketBig{ \alpha , U^n \beta  }_{ L^2(\mu) } = 
		\begin{cases}
			0 & \mbox{ if } \alpha \in \Disc, \beta \in \Disc^\bot \\
			O(1) & \mbox{ if } \alpha, \beta \in \Disc \\
			\mbox{Cesaro summable} & \mbox{ if }  \alpha, \beta \in \Disc^\bot 
		\end{cases}
	\end{equation}
	The last relation can be stated more accurately as
	\begin{equation} \label{eqn:untry:3}
		\lim_{N\to \infty} \frac{1}{N} \sum_{n=0}^{N-1} \bracketBig{ \alpha , U^n \beta  }_{ L^2(\mu) } = 0 .
	\end{equation}
	This phenomenon is known as \emph{weak decay of correlations}. If the component $\Disc$ is trivial then the system is called \emph{weak mixing} \cite[e.g.]{DasGiannakis_delay_2019}. For any two observables $\alpha, \beta$ drawn from $\Disc^\bot$, their inner products converge to $0$ in a sense made precise in \eqref{eqn:untry:3}. This plays an important role in any analysis that involves projections to subspaces. 
	
	\paragraph{Linear hypothesis space} Any numerical method attempting to determine the feedback function $w_k$ must chose a set of candidate functions. This set is known as the \emph{hypothesis space}. In a parametric approach, this set is a parametric family of functions, corresponding to parameter values chosen from over some range of values. In a data-driven approach the hypothesis space $\mathcal{H}$ is usually a finite dimensional linear space, spanned by a basis $h_1, \ldots, h_m$. In that case the collection
	\begin{equation} \label{eqn:def:W}
		\mathcal{W} := \spn\SetDef{ h_i \circ \Phi_l }{ 1\leq i\leq m, \, \leq l\leq L }
	\end{equation}
	is a finite subspace of $L^2(\mu)$, and the estimate $\hat{w}_k$ satisfies
	\begin{equation} \label{eqn:kcnv83}
		\hat{w}_k \circ \Phi = \proj_{\mathcal{W}} U^k \phi .
	\end{equation}
	For example, if the hypothesis space is restricted to space $\Linear\left( \real^L; \real^d \right)$ of linear functions, then $\mathcal{W} = \spn \Phi$. In the rest of this paper, we shall focus on this scenario where the hypothesis space is linear. We state this formally :
	
	\begin{Assumption} \label{A:hypoth}
		The hypothesis space $\mathcal{W}$ is a finite dimensional subspace of $L^2(\mu)$, and contains the constant function $1_{\real^L}$.
	\end{Assumption}
	
	In most learning techniques, a bias or offset constant is calculated separately, thus satisfying the criterion that $\mathcal{W}$ contains constant functions.
	
	\paragraph{Error from direct forecasts} Let $\pi$ denote the projection $\proj_{\mathcal{W}}$, and set $\Delta := \Id-\pi$. For ease of notation, we will denote $\hat{w}_1$ simply by $\hat{w}$, in the rest of this section. Define the \emph{projection error}  to be the quantity
	\begin{equation} \label{eqn:def:delta}
		\delta = \delta(\mathcal{H}) := \norm{ \Delta U\phi }_{L^2(\mu)} .
	\end{equation}
	This is the component of the measurement $\phi$ not recoverable using our choice of hypothesis space. Note that as the size of the hypothesis space increases, $\delta$ converges to $0$. 
	We shall first examine the performance of the direct forecast method. 
	
	\begin{theorem}[Error from direct forecast] \label{thm:direct}
		\cite[Thm 4.1]{BerryDas_learning_2022} Let Assumptions~\ref{A:f}, \ref{A:data} and \ref{A:pPhi} hold, and assume the notations in \eqref{eqn:def:hat_w}, \eqref{eqn:def:feedback_1} and \eqref{eqn:feedback_1_init}. Let $\delta$ be as in \eqref{eqn:def:delta}. Then the error from direct iteration is given by
		\begin{equation} \label{eqn:direct:1}
			\text{error}_{\text{direct}}(n) = \norm{ \left( \Id - \pi \right) U^n\phi }_{L^2(\mu)} .
		\end{equation}
		Now assume that Assumption~\ref{A:hypoth} holds. Then there is a subset $\num'\subseteq \num$ with density $1$ such that the following holds : 
		\begin{enumerate} [(i)]
			\item For every $\epsilon>0$, if the hypothesis space $\mathcal{W}$ is chosen large enough, then
			\[ \lim_{n\in \num', n\to\infty} \text{error}_{\text{direct}}(n) = \norm{ \phi - \proj_{\Disc} \phi }_{L^2(\mu)} + \epsilon. \]
			\item If $f$ is weakly mixing, then for every choice of $\mathcal{W}$
			\[ \lim_{n\in \num', n\to\infty} \text{error}_{\text{direct}}(n) = \text{var}_{\mu} := \norm{ \phi - \mu(\phi)}_{L^2(\mu)} . \]
			\item If $f$ is strongly mixing, the set $\num'$ can be taken to be the entire set $\num$.
			\item If $f$ has purely discrete spectrum, then for every $\epsilon>0$, if the hypothesis space $\mathcal{W}$ is chosen large enough, then
			\[ \text{error}_{\text{direct}}(n) < \epsilon, \quad \forall n\in\num . \]
		\end{enumerate}
	\end{theorem}
	
	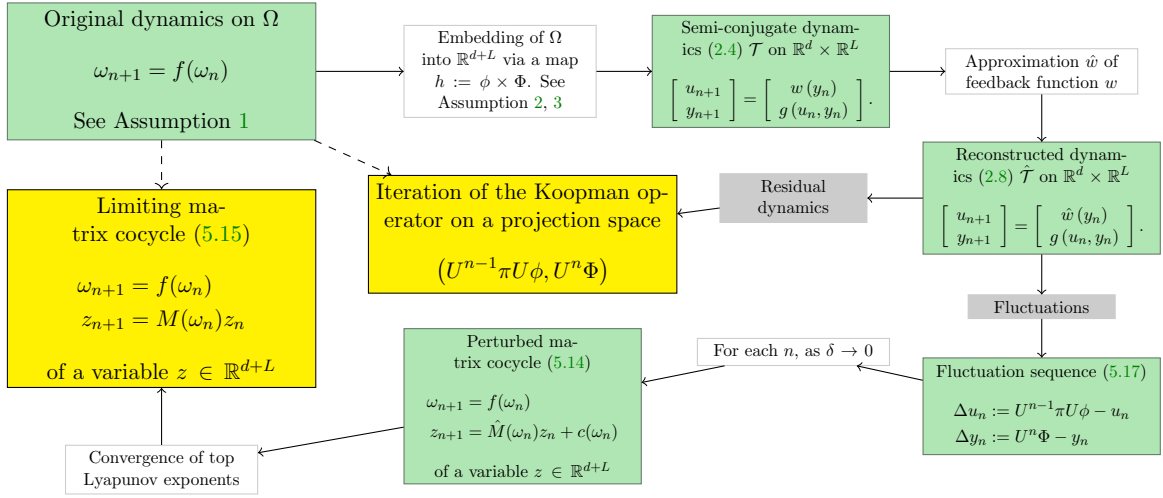
\begin{figure}\center
		\begin{tikzpicture}[scale=0.6, transform shape]
			\node [style={rect11}, scale=1.3] (N1) at (-0.3\columnA, 0.4\rowA) { Original dynamics on $\Omega$ \[ \omega_{n+1} = f(\omega_n) \] See Assumption \ref{A:f} } ;
			\node [style={rect2}] (N2) at (1.2\columnA, 0.4\rowA) { Embedding of $\Omega$ into $\real^{d+L}$ via a map $h:= \phi\times \Phi$. See Assumption \ref{A:data}, \ref{A:pPhi} } ;
			\node [style={rect11}] (N3) at (2.4\columnA, 0.4\rowA) { Semi-conjugate dynamics \eqref{eqn:def:feedback_1} $\reconstruct $ on $\real^d \times \real^L$ 
				\[\left[ \begin{array}{c} u_{n+1} \\ y_{n+1} \end{array} \right]
				= \left[ \begin{array}{c} w \left( y_n \right) \\ g\left( u_{n}, y_n \right) \end{array} \right].\] 
			} ;
			\node [style={rect2}] (N4) at (3.6\columnA, 0.4\rowA) { Approximation $\hat{w}$ of feedback function $w$ } ;
			\node [style={rect11}] (N5) at (3.6\columnA, -1\rowA) { Reconstructed dynamics \eqref{eqn:def:feedback_2} $\hat{\reconstruct} $ on $\real^d \times \real^L$ 
				\[\left[ \begin{array}{c} u_{n+1} \\ y_{n+1} \end{array} \right]
				= \left[ \begin{array}{c} \hat{w} \left( y_n \right) \\ g\left( u_{n}, y_n \right) \end{array} \right].\] 
			} ; 
			\node [style={rect4}] (N6) at (2.5\columnA, -1\rowA) { Residual dynamics } ;
			\node [style={rect12}, scale=1.3] (N7) at (1.3\columnA, -1.4\rowA) { Iteration of the Koopman operator on a projection space \[ \paran{ U^{n-1} \pi U \phi , U^n \Phi } \] } ;
			\node [style={rect4}] (N8) at (3.6\columnA, -2.2\rowA) { Fluctuations } ; 
			\node [style={rect11}] (N9) at (3.6\columnA, -3.3\rowA) { Fluctuation sequence \eqref{eqn:def:Delta_uy} 
				\begin{align*}
					\Delta u_n &:= U^{n-1} \pi U \phi - u_n \\
					\Delta y_n &:= U^n\Phi - y_n
				\end{align*} 
			} ;
			\node [style={rect2}] (N10) at (2.5\columnA, -2.7\rowA) { For each $n$, as $\delta \to 0$ } ; 
			\node [style={rect11}] (N11) at (1.3\columnA, -3.3\rowA) { Perturbed matrix cocycle \eqref{eqn:ab_approx}
				\begin{align*}
					\omega_{n+1} &= f(\omega_n) \\
					z_{n+1} &= \hat{M}( \omega_n ) z_n + c(\omega_n)
				\end{align*} 
				of a variable $z \in \real^{d+L}$} ;
			\node [style={rect2}] (N12) at (-0.3\columnA, -4\rowA) { Convergence of top Lyapunov exponents } ; 
			\node [style={rect12}, scale=1.3] (N13) at (-0.3\columnA, -2\rowA) { Limiting matrix cocycle \eqref{eqn:ab_exact}  
				\begin{align*}
					\omega_{n+1} &= f(\omega_n) \\
					z_{n+1} &= M( \omega_n ) z_n 
				\end{align*} 
				of a variable $z \in \real^{d+L}$} ;
			\draw[-to] (N1) to (N2);
			\draw[-to] (N2) to (N3);
			\draw[-to] (N3) to (N4);
			\draw[-to] (N4) to (N5);
			\draw[-to] (N5) to (N6);
			\draw[-to] (N6) to (N7);
			\draw[-to] (N5) to (N8);
			\draw[-to] (N8) to (N9);
			\draw[-to] (N9) to (N10);
			\draw[-to] (N10) to (N11);
			\draw[-to] (N11) to (N12);
			\draw[-to] (N12) to (N13);
			\draw[dashed,->] (N1) -- (N7);
			\draw[dashed,->] (N1) -- (N13);
		\end{tikzpicture}
		\caption{Iteration error. The flowchart presents an analysis of the error resulting in an iterative forecasting of a dynamical system. The green boxes represent various conceptual dynamical systems needed to interpret the result. The white boxes mentions in brief the theoretical connections between these links. The target is a dynamical system shown in the top left corner. The main realization of Section \ref{sec:predict} and Theorem \ref{thm:iterative} is that the performance of the prediction can be expressed as the sum of two separate parts, as shown in yellow boxes above. The residual components can be described using purely by operator theoretic means, as an iteration of the Koopman operator $U$ on a subspace. The other component represents the fluctuations from this simple form. These fluctuations can be modeled by a special dynamical system called a perturbed matrix cocycle. The residual component can be studied well from the operator theoretic properties of  $U$. The growth of the fluctuations can be understood well from the ergodic theory of matrix cocycles.}
		\label{fig:iter_error}
	\end{figure} 
	
	\paragraph{Role of the spectrum} A deeper look into dynamical systems theory thus reveals how secondary characteristics such as the Koopman operator, and even tertiary characterics such as $\Disc, \Disc^\bot$ have a direct bearing on the performance of numerical methods. The error of direct forecasts increase at the rate of mixing of the system. If there are quasiperiodic components in the dynamics, the direct method is more effective in retaining that component. The term spectral analysis is commonly used in signal processing. One can come to similar conclusions as signal processing based on an ergodig interpretation of spectrum. 
	Given any nonzero function $\psi \in L^2(\mu)$, we define its normalized autocorrelation (w.r.t. the underlying dynamics) as
	\[ \AutCor(n; \psi):= \norm{\psi}^{-2} \left\langle U^n \psi, \psi \right\rangle . \]
	Autocorrelation is a statistical property of signals used frequently in classical timeseries analysis \citep[e.g.][]{box2015time}. It can be shown that \cite[rem 4.5]{BerryDas_learning_2022}
	\begin{equation} \label{eqn:phi_hypo_autocorr}
		\phi\in\mathcal{W} \imply \text{error}_{\text{direct}}(n)^2 \leq \norm{\phi}^2 \left[ 1 - \AutCor(n; \phi)^2 \right].
	\end{equation}
	Thus, if the hypothesis space happens to include the initial observation map $\phi$, then the growth of the direct error is directly related to the autocorrelation function of the observed signal $\phi$. Equation~\ref{eqn:phi_hypo_autocorr} thus combines concepts from learning theory, ergodic theory and timeseries analysis.
	We now turn our attention towards analyzing the error from iterative prediction. We begin by reviewing the concept of matrix cocycles.
	
	\paragraph{Matrix cocycles} Let Assumption~\ref{A:f} hold, and $G : \Omega\to GL(\real, m)$ be a measurable map. Then it generates a \emph{matrix cocycle} (see \cite{FroylandEtAl_coherent_2010}, \citep[][3.4]{Arnold_random_1991}), which is the map
	\begin{equation}\label{eqn:ih7b}
		\cocyc:\Omega\times \num_0\to GL(\real;m), \quad \cocyc(n, \omega) := 
		\begin{cases}
			\Id_d &\mbox{if } n=0\\
			G(f^{n-1} \omega ) \cdots G(\omega) &\mbox{if } n>0 \\
			G( f^{-|n|} \omega)^{-1} \cdots G(f^{-1} \omega )^{-1} &\mbox{if } n<0
		\end{cases} .
	\end{equation}
	$\cocyc$ is called a $GL(m;\real)$-valued cocycle over the dynamics $(\Omega,\mu,f)$ generated by $f$. It has the property
	\begin{equation}\label{eqn:def:cocyc}
		\cocyc(m+n, \omega) = \cocyc(n, f^{m} \omega) \cdot \cocyc(m, \omega), \quad \forall \omega\in \Omega, \quad \forall m,n\in \integer. 
	\end{equation}
	Here the $\cdot$ notation denotes the matrix multiplication. Equation \eqref{eqn:def:cocyc} is the defining equation of a matrix cocycle. Conversely, given any map $\mathcal{G}:\Omega\times\num_0\to GL(m;\real)$ satisfying \eqref{eqn:def:cocyc}, one has a generator $G : \Omega \to GL(m;\real)$ so that $\cocyc$ is related to $G$ via \eqref{eqn:ih7b}. One of the immediate consequences of \eqref{eqn:def:cocyc} is that
	\[ \cocyc(0, \omega) = \Id_m, \quad \cocyc( -n, \omega ) = \cocyc( n, f^{-n} \omega )^{-1}, \quad \forall \omega\in \Omega . \]
	When the initial point $\omega_0\in\Omega$ is fixed, we will drop it from the notation and define
	\[ \cocyc(n-1,j) := \cocyc(n-j, f^{j} \omega_0) = G(f^{n-1} \omega_0) \cdots G(f^{j} \omega_0) . \]
	Matrix valued cocycles arise naturally in multiple ways in dynamical systems. For example, if $\Omega$ is a $m$-dimensional manifold and $f$ a differentiable map, then $\cocyc( \omega,n) := Df^n(\omega)$ is a $GL(m;\real)$ cocycle. It tracks the rate of growth of deviations from a reference trajectory. Matrix cocycles arise in our study from the following matrix valued functions
	\begin{equation}\begin{split} \label{eqn:def:WG_matrices}
			W :\Omega\to \real^{d\times L}, & \quad W(\omega) := Dw \circ \Phi(\omega) ,\\
			\hat{W} :\Omega\to \real^{d\times L}, & \quad \hat{W} (\omega) := D\hat{w} \rvert_{ \Phi(\omega)} = D\hat{w} \circ \Phi(\omega),\\
			G^{(1)} :\Omega\to \real^{L\times d}, & \quad G^{(1)}(\omega) := \nabla_1 g \rvert_{h(\omega)} = \nabla_1 g \circ h(\omega),\\ 
			G^{(2)} :\Omega\to \real^{L\times L} , & \quad G^{(2)}(\omega) := \nabla_2 g \rvert_{h(\omega)} \nabla_2 g \circ h(\omega),
	\end{split}\end{equation}
	The matrix valued functions $W$ and $\hat{W}$ stem from the true and approximated feedback functions $w$ and $\hat{w}$ respectively. These combine to give two matrix valuef functions $M, \hat{M} : \Omega\to \real^{(L+d)\times (L+d)}$
	and their combination 
	\[\hat{M}(\omega) := \left[\begin{array}{cc} 0^{d\times d} & \hat{W}(\omega) \\ G^{(1)}(\omega) & G^{(2)}(\omega)  \end{array}\right] , \quad 
	M(\omega) := \left[\begin{array}{cc} 0^{d\times d} & W(\omega) \\ G^{(1)}(\omega) & G^{(2)}(\omega)  \end{array}\right]\]
	Next consider the vector-valued functions
	\[ c:\Omega\to \real^{d+L}, \, c(\omega) := \paran{ \vec{0}_d \,,\, G^{(1)}(\omega) \left( U^{-1} \Delta \phi \right)(\omega) }. \]
	These various matrix and vector valued functions can be used to build a non-autonomous dynamical system on $\real^{d+L}$. Fix an $\omega\in \Omega$ and define
	\begin{equation}  \label{eqn:ab_approx}
		\left[\begin{array}{c} a_{n+1} \\ b_{n+1} \end{array} \right] = \hat{M} \left( f^n \omega \right) \left[\begin{array}{c} a_{n} \\ b_{n} \end{array} \right] +c\left( f^n\omega \right) , \quad a_1=0^d, \; b_0=0^L .
	\end{equation}
	We call such a system a \emph{perturbed matrix cocycle}. Note that as the size of the hypothesis space is increased, the function $c$ converges to $0$ in $L^2(\mu)$ norm, and the dynamics of $(a_n, b_n)$ gets closer to that of the matrix cocycle generated by $M$. Thus the limiting dynamics is
	\begin{equation}  \label{eqn:ab_exact}
		\left[\begin{array}{c} a_{n+1} \\ b_{n+1} \end{array} \right] = M \left( f^n \omega \right) \left[\begin{array}{c} a_{n} \\ b_{n} \end{array} \right]
	\end{equation}
	The dynamics in both \eqref{eqn:ab_approx} and \eqref{eqn:ab_exact} depends on the initial state $\omega$. If $\omega$ is allowed to vary, then $a_n, b_n$ become functions of $\omega$. We shall overuse notation and also denote these functions as $a_n, b_n$. 
	
	\paragraph{Growth of the iterative error} Figure \ref{fig:iter_error} presents an overview of the analysis of the iterative error. The targeted dynamics is the system $(\Omega, f, \mu)$ as described in Assumption Section \ref{sec:learn} presents a common framework for all learning techniques for the dynamics law $f$. In a data-driven approach the original phase space $\Omega$ is absent in the calculation. Instead one tries to reconstruct the dynamics in $\real^{d+L}$, through the map $\reconstruct$ from \eqref{eqn:def:feedback_1}. The dynamics is learnt by approximating the feedback function $w$ as some function $\hat{w}$. The result is an approximated dynamics map $\hat{\reconstruct}$ as in \eqref{eqn:def:feedback_2}. The study of the iteration error is thus the deviation of the first $d$ coordinates of the dynamics under $\hat{\reconstruct}$ from that of $\reconstruct$. We wish to make comparison
	\begin{equation} \label{eqn:pd0g}
		\left[\begin{array}{c} u_n \\ y_n \end{array} \right] := \hat{\reconstruct}^n \left( \begin{array}{c} \phi \\ \Phi \end{array} \right) \quad \mbox{v.s.} \quad \left[\begin{array}{c} U^{n-1} \pi U \phi \\ U^n\Phi \end{array} \right]
	\end{equation}
	In other words, we compare the iterations of \eqref{eqn:def:feedback_2} with iterations of $U$ on a pair of functions. The comparison \eqref{eqn:pd0g} has two components to it :
	
	\begin{enumerate}
		\item \textbf{Reference dynamics.} The sequence generated by $\paran{U^{n-1} \pi U \phi, U^n\Phi}$ is merely the evolution of the functions $\pi \phi$ and $\Phi$ under the action of $U$. This sequence remains bounded in $L^2(\mu)$ sense. The simplicity of this sequence makes it a good reference to compare the performance of the forecast against. Note that the first $d$ coordinates of the latter sequence differs from the true sequence $\hat{\omega}_n$ \eqref{eqn:def:hatomega} via the map 
		\[ U^n \circ \paran{\Id - \pi} \phi .\]
		The $L^2(\mu)$ norm of this function is constant and equal to that of $\Delta \phi = \paran{\Id - \pi} \phi$. This function $\Delta \phi$ represents the residual error created during the choice of the hypothesis space $\calH$. If $\calH$ is made larger, the residual error also decreases. The residual error is something that is determined from the choice of hypothesis space. The choice of the feedback function cannot decrease this error.
		\item \textbf{Fluctuations.} The comparison in \eqref{eqn:pd0g} leads to the study of their differences, which form the sequence 
		\begin{equation} \label{eqn:def:Delta_uy}
			\left[\begin{array}{c} \Delta u_n \\ \Delta y_n \end{array} \right] = \left[\begin{array}{c} U^{n-1} \pi U \phi \\ U^n\Phi \end{array} \right] - \left[\begin{array}{c} u_n \\ y_n \end{array} \right], \quad \forall n\in\num_0.
		\end{equation}
		We call this sequence the fluctuations from the reference system.
	\end{enumerate}
	
	Equation \eqref{eqn:pd0g} suggests that a good reference or comparison for the iterative error is the Koopman dynamics. To establish this concretely one needs to bound the growth of the fluctuations. Their dynamics can be bounded by describing them as a perturbed matrix cocycle. More precisely : 
	
	\begin{theorem}[Error from iterative forecast] \label{thm:iterative}
		\cite[Thm 4.2]{BerryDas_learning_2022} Let Assumptions~\ref{A:f}, \ref{A:data} and \ref{A:pPhi} hold, and assume the notations in \eqref{eqn:def:hat_w}, \eqref{eqn:def:feedback_1} and \eqref{eqn:feedback_1_init}. Fix an initial state $\omega\in \Omega$ and let $(u_n, y_n)$ be successive iterations of the system \eqref{eqn:def:feedback_1}, while $(a_n, b_n)$ be iterations of the dynamics in \eqref{eqn:ab_approx}. Let $\delta$ be as in \eqref{eqn:def:delta}.
		\begin{enumerate} [(i)]
			\item  The deviations \eqref{eqn:def:Delta_uy} are related to the iterations \eqref{eqn:ab_approx} of the perturbed matric cocycle as :
			\begin{equation} \label{eqn:u_vs_a}
				\Delta u_n = a_n + \bigO{a_{n-1}}^2, \quad \lim_{\delta\to 0} \frac{ \norm{\Delta u_n} }{ \norm{a_n} } = 1 .
			\end{equation}
			\item Suppose that dynamics $(\Omega, \mu, f)$ has the property of $L^2$ Pesin sets. Let $\lambda_1 = \lambda_1(\mathcal{M})$  the maximal Lyapunov exponent of the cocycle generated by $\hat{M}$. Then for every $\epsilon>0$, there is a a constant $C_{\epsilon}$ such that
			\begin{equation} \label{eqn:def:iter_bound_L2}
				\text{error}_{\text{iter}}(n) = \norm{\Delta u_n}_{L^2(\mu)} = \delta C^{(2)}_{\epsilon} \bigO{ e^{n (\lambda_1+\epsilon)} }, \quad \mbox{as } n\to\infty .
			\end{equation}
			for a constant $C_{\epsilon}$ that depends only on $\epsilon$.
		\end{enumerate}
	\end{theorem}
	
	The idea of Pesin sets and their regularity are ways to characterize the rate of non-uniformity in in chaotic systems. See \cite[Sec 11]{BerryDas_learning_2022} for more details. Thus the fluctuations grow at a rate, whose asymptotic exponential rate of growth is related to the Lyapunov exponent of the matrix cocycle generated by $M$ \eqref{eqn:ab_exact}. 
	
	\paragraph{Direct vs Iterative} The explicit formulas for the direct and iterative schemes reveal a basic mathematical law that makes the direct method unsuitable for long term prediction. It involves the operator $\pi U^n$, which projects the evolving measurement $\phi$ back into the space $\mathcal{W}$. For strongly mixing systems, the Koopman operator drives out any function from any finite dimensional subspace, up to a constant function. On the other hand, the iterative method always involves the term $U^{n-1} \pi U$.  The crucial difference is that the projection $\pi$ is not made after the application of $U^n$, but always to the static operator $U$. The $U^{n-1}$ in front of the $\pi$ then merely acts as a rotation / unitary transform and thus does not change the $L^2(\mu)$ (i.e., RMS) magnitude of the error. 
	
	\paragraph{Overfitting error vs projection error} Equation~\eqref{eqn:def:iter_bound_L2} shows that the projection rate grows exponentially as expected from the presence of a Lyapunov exponent. The rate of growth is proportional to the smoothness of the learnt function $\hat{w}$, while the multiplicative constant is proportional to the projection error $\delta$. Thus this displays a trade-off between projection error and overfitting, one can minimize the projection error by increasing the hypothesis space. But the resulting learnt function may be too oscillatory, as a result increasing the instability of the feedback system \eqref{eqn:def:feedback_1}. On the other hand, if one approximates $w_1$ by a less oscillatory function, our base error $\delta$ itself will be large to begin with. To state this trade-off more precisely, define
	\[ \theta(\epsilon) := \inf \SetDef{ \norm{D\hat{w}} }{ \hat{w} \in \mbox{ some hypothesis space } \mathcal{H}, \, \delta(\mathcal{H}) < \epsilon  } .\]
	Then $\theta$ is a non-decreasing function of $\epsilon$, satisfying
	\[ \theta\left( \norm{\phi} \right) = 0, \quad \lim_{\epsilon\to 0^+} \theta(\epsilon) = \norm{Dw}. \]
	Thus \eqref{eqn:def:iter_bound_L2} can be rewritten as 
	\[ \text{error}_{\text{iter}}(k) = \epsilon \bigO{ k \theta(\epsilon)^k }, \quad \mbox{as } k\to\infty. \]

	This completes the statement of our main results. We next look at some numerical experiments to verify the statements of Theorems \ref{thm:iterative} and \ref{thm:direct}. 
	
	\section{Numerical experiments} \label{sec:expt}
	
	\begin{figure}[!ht]\center
		\includegraphics[width=.48\linewidth]{\figs 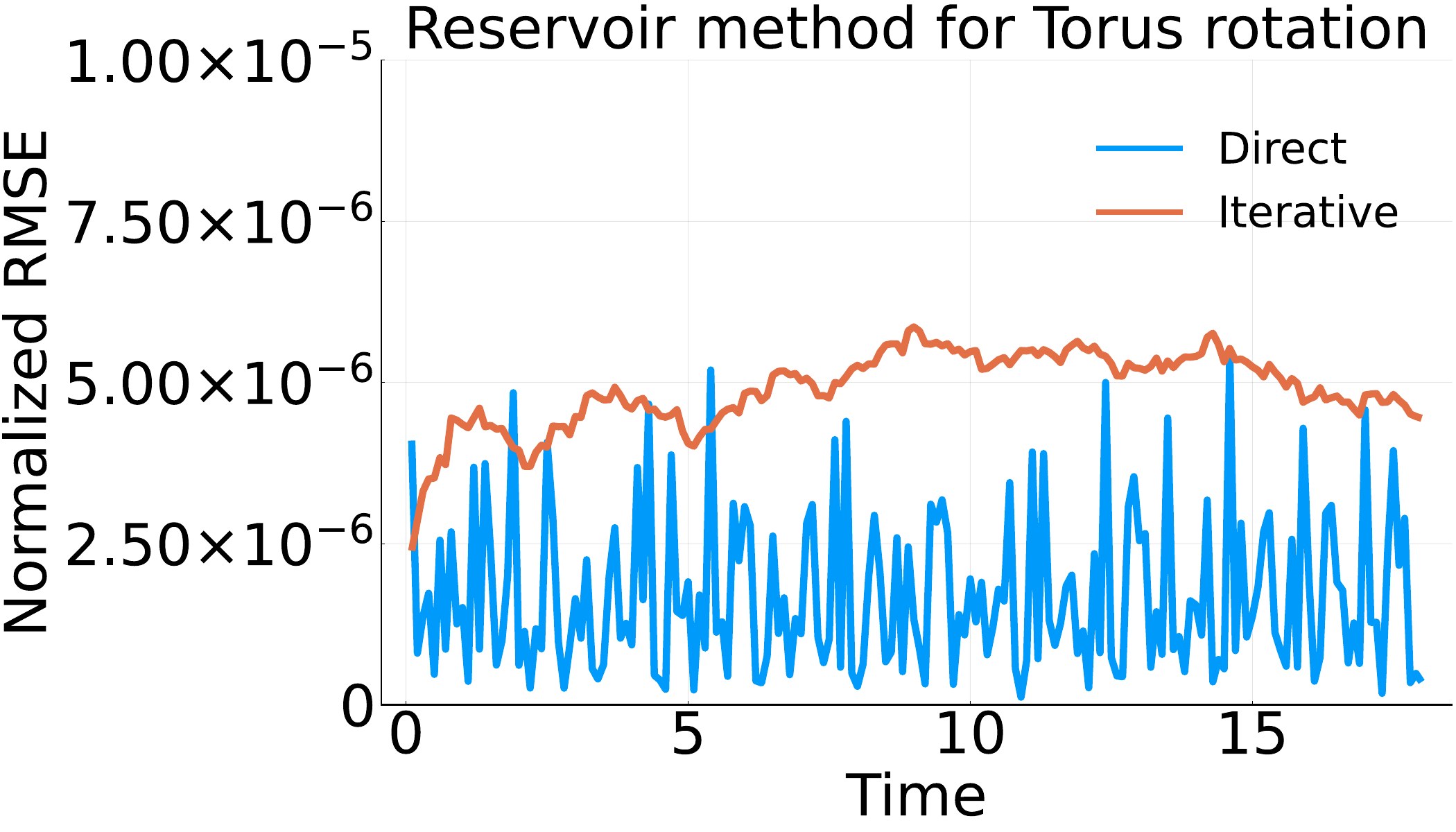}
		\includegraphics[width=.48\linewidth]{\figs 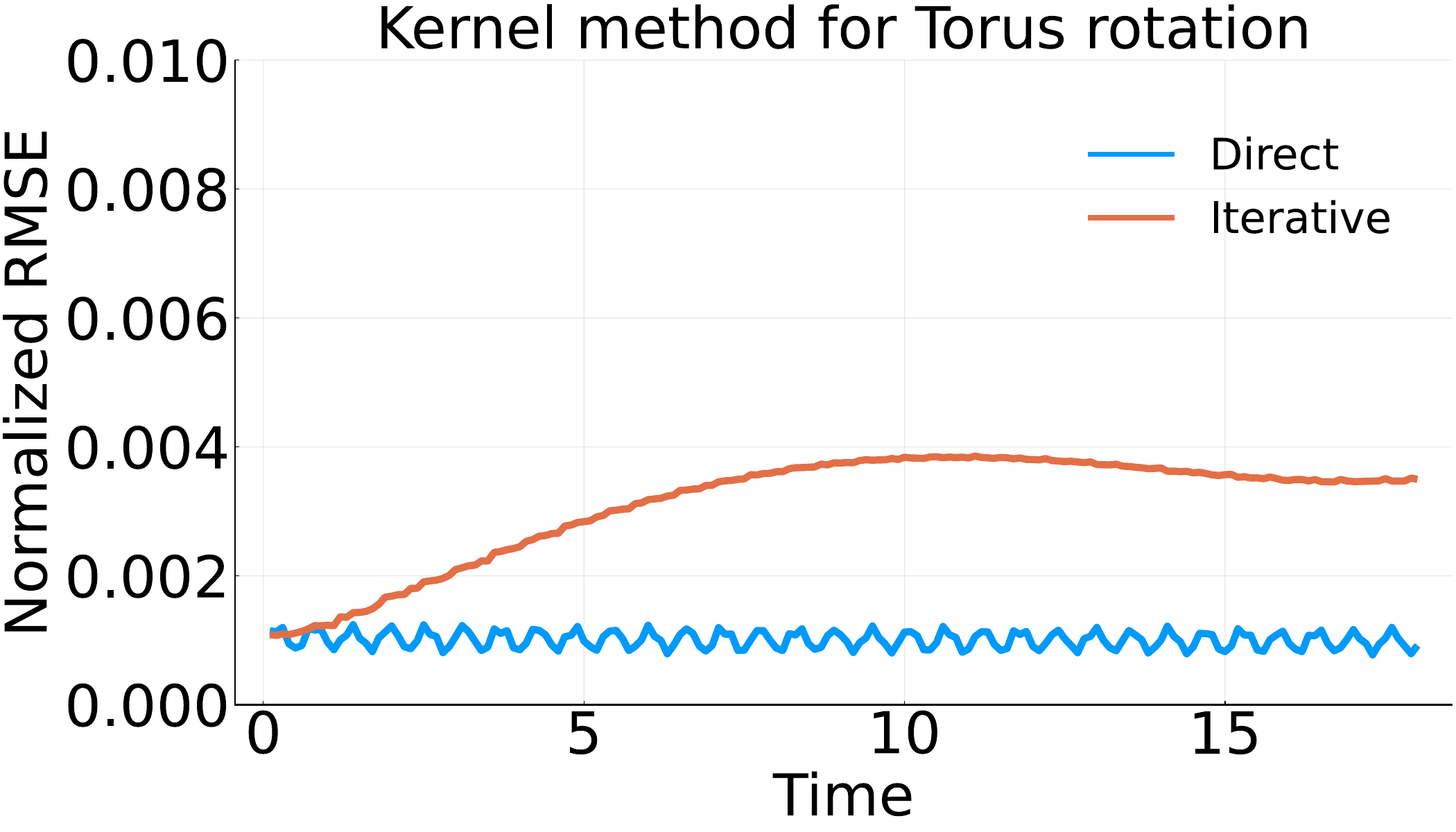}
		\includegraphics[width=.48\linewidth]{\figs 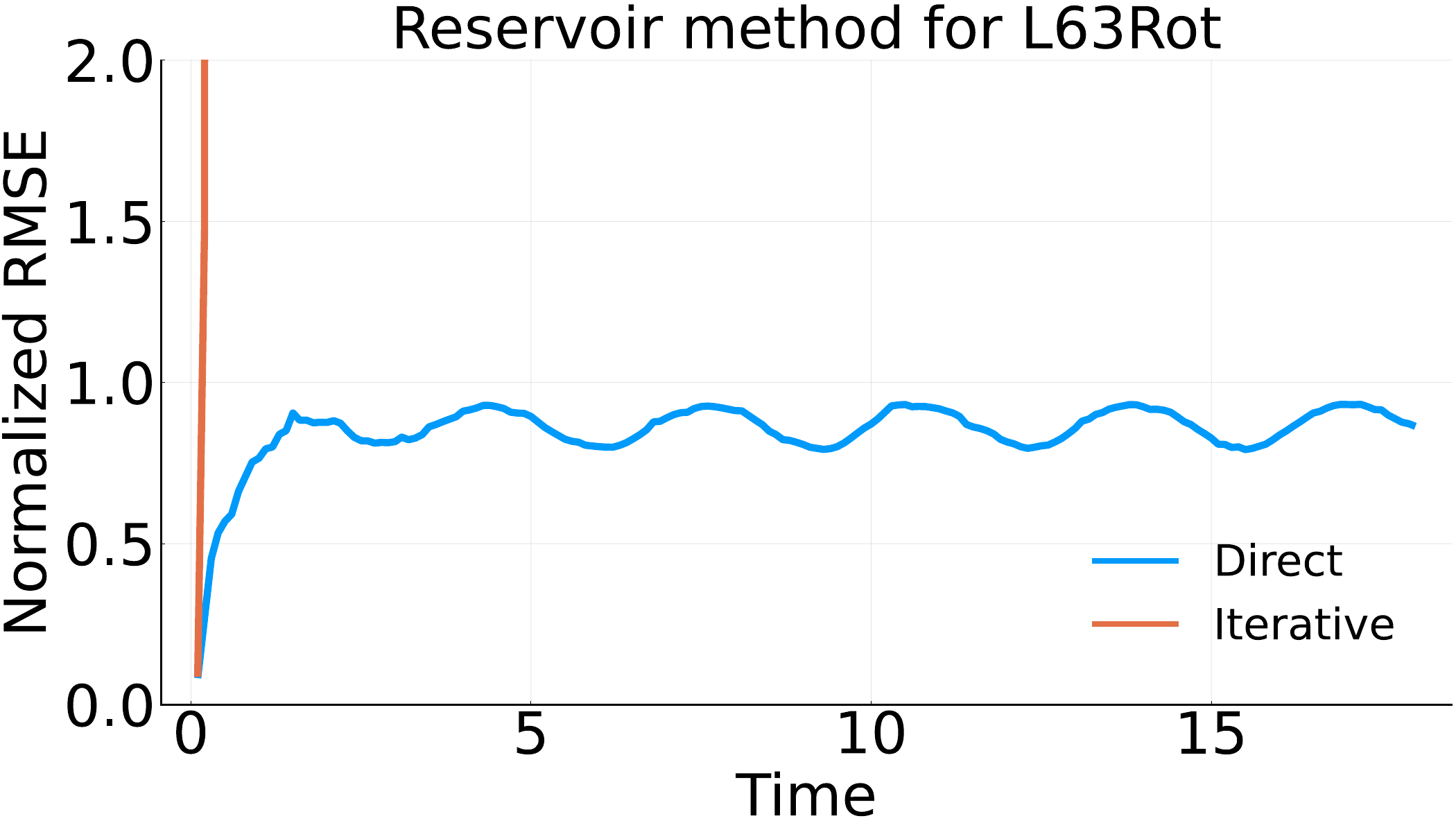}
		\includegraphics[width=.48\linewidth]{\figs 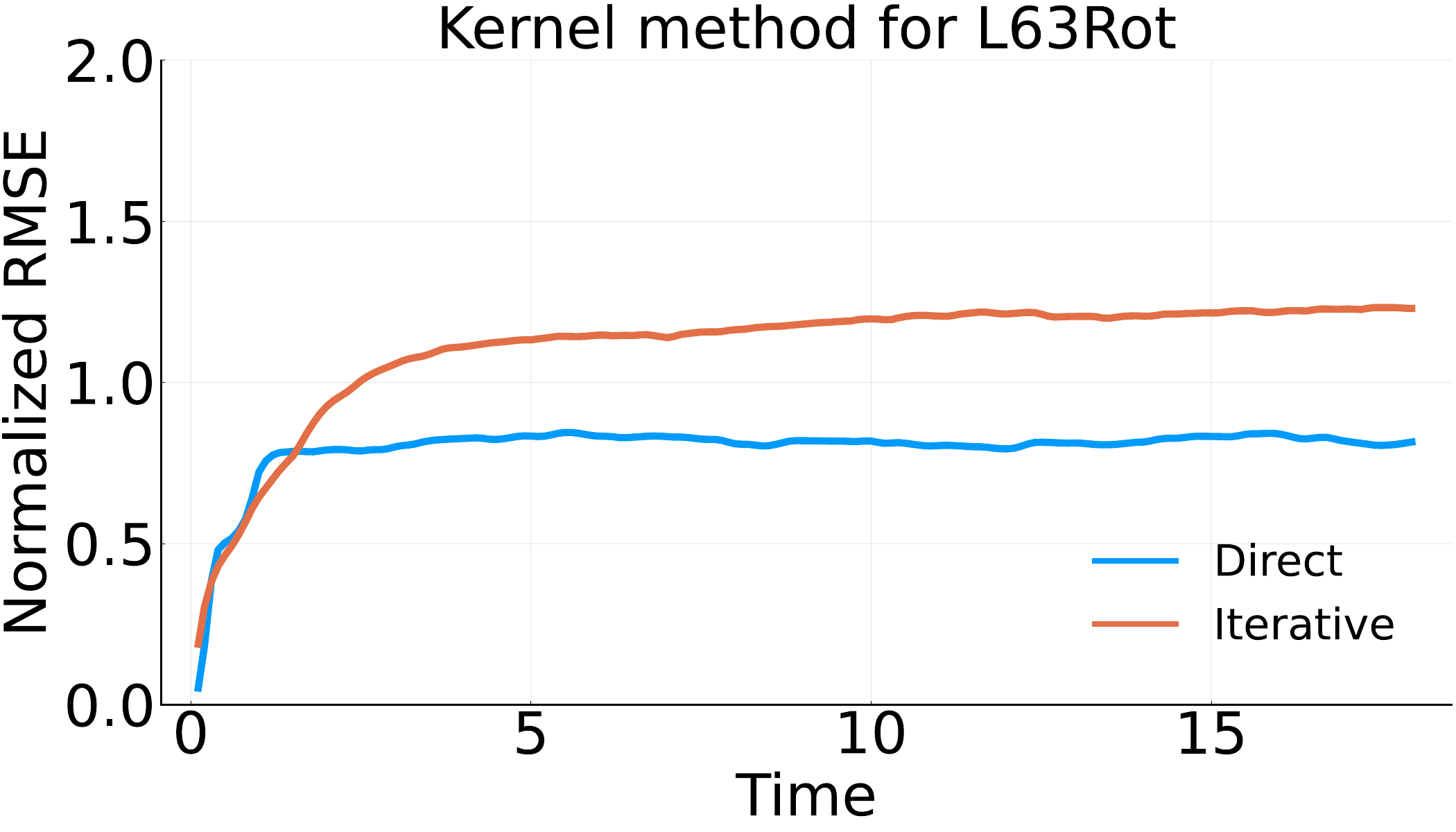}
		\caption{Performance of the two reconstruction techniques for (i) a quasiperiodic rotation on the torus $\mathbb{T}^2$ (bottom panels); and (ii) a Cartesian product of L63 with a simple harmonic oscillator (top panels). By Theorem~\ref{thm:direct}~(iv), if one has a proper embedding and a good approximation $\hat{w}$ of $w$ one can achieve arbitrarily small errors for the torus rotation,  for all forecast times. This is supported by the fact that the direct methods for both the paradigms show errors of the order of $10^{-6}$. Since the torus rotation has all Lyapunov exponents zero, by Theorem~\ref{thm:iterative}~(ii) and Theorem~\ref{thm:lambda1}~(ii), the error from the iterative techniques should grow sub-exponentially, as supported by the figures. The system \eqref{eqn:L63Rot} is a mixed spectrum system, i.e., the splitting in \eqref{eqn:def:L2split} is non-trivial.  }
		\label{fig:T2_L63Rot}
	\end{figure}
	
	\begin{figure}\center
		\includegraphics[width=.48\linewidth]{\figs 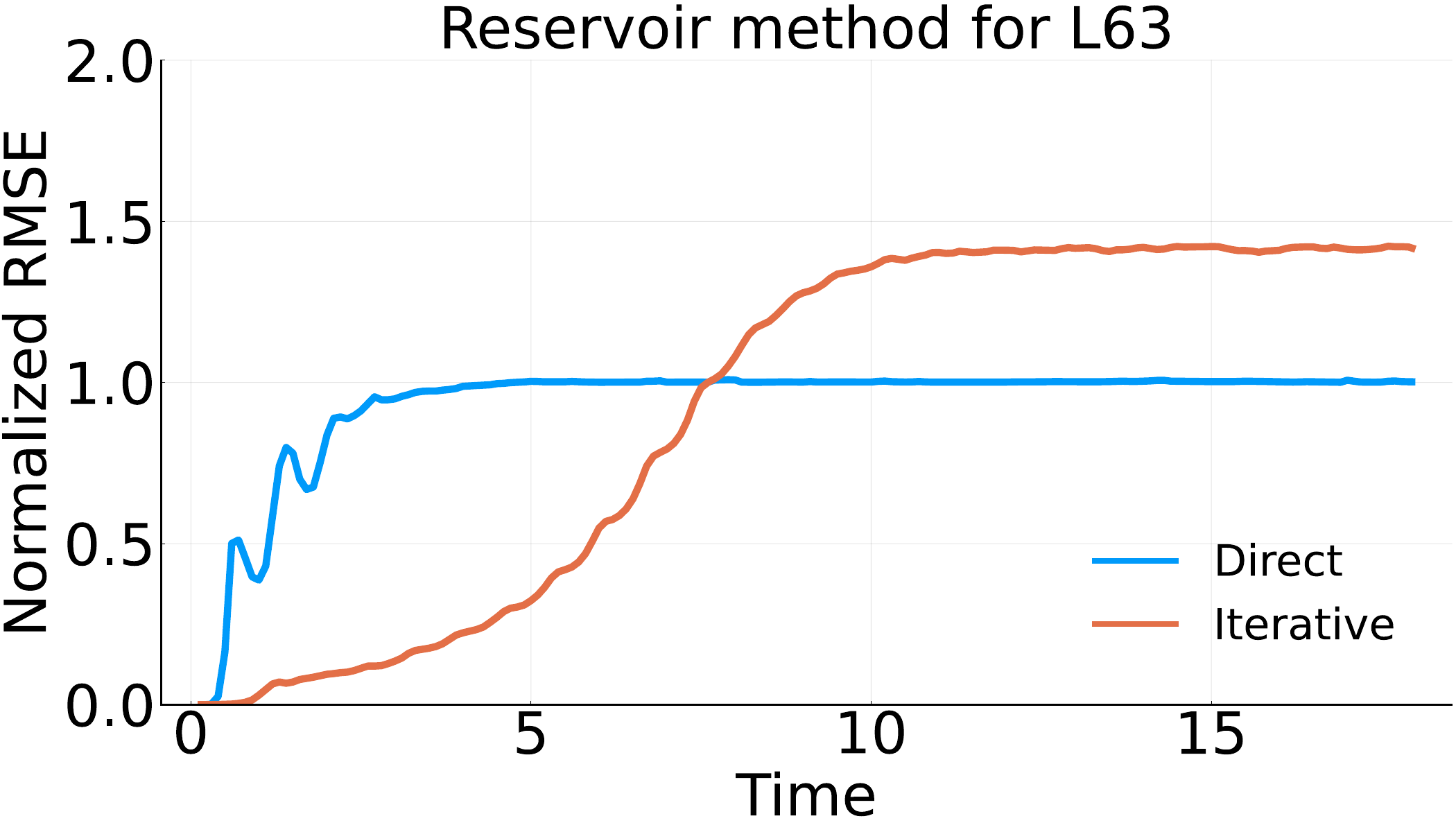}
		\includegraphics[width=.48\linewidth]{\figs 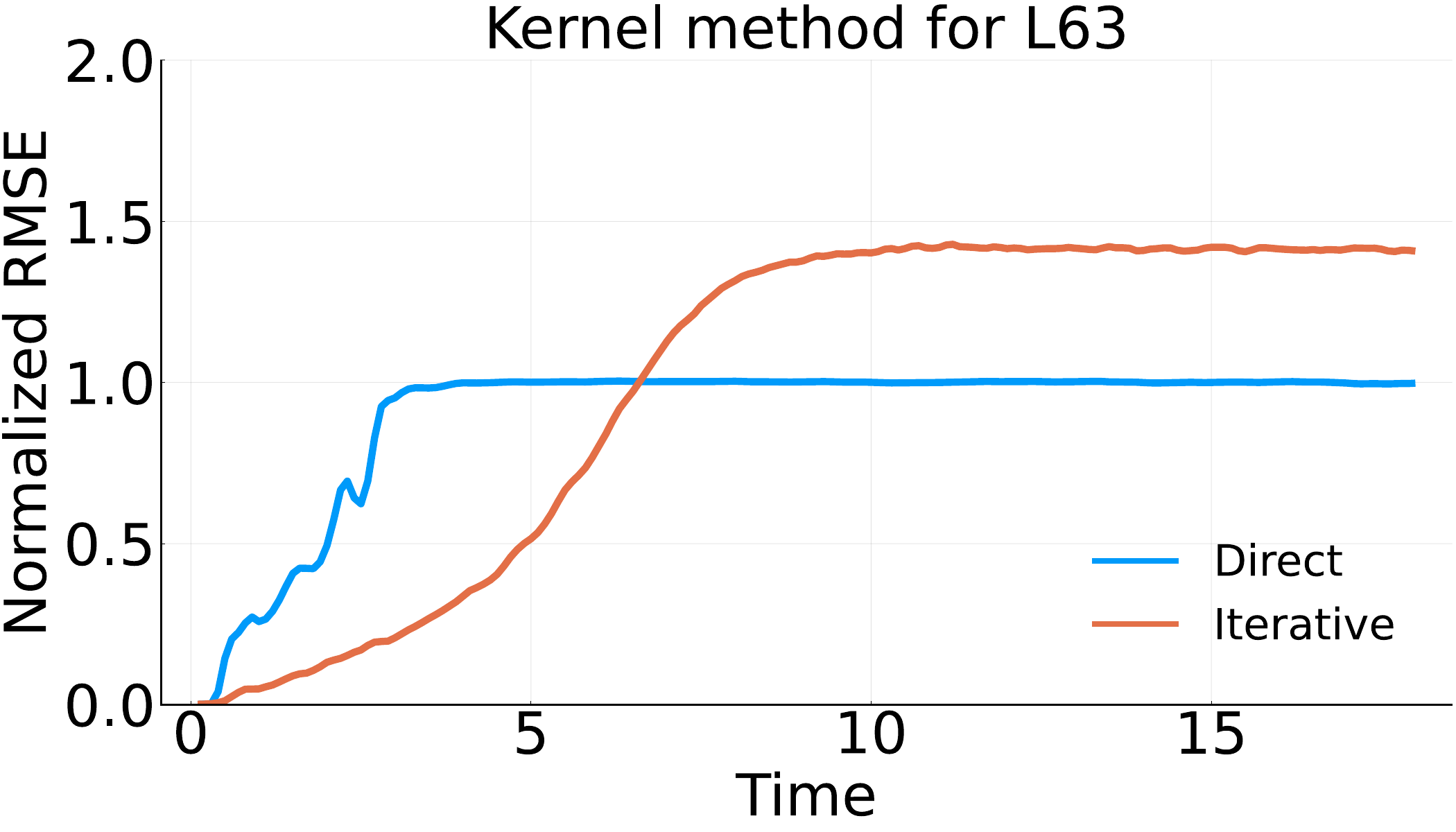}
		\caption{Results of forecasting attempts on the Lorenz63 system, using the invariant-graph paradigm on the left, and the delay-coordinate based paradigm on the right. For both these paradigms, we compare the performance of direct and iterative modes of forecasting (see \eqref{eqn:def:err_iter}, \eqref{eqn:def:err_direct} for definitions). The horizontal axis shows the forecast-time $n$, and the vertical axis is the RMS error of forecast as a function of $n$, for a signal of unit $L^2(\mu)$-norm. The RMS error is meant to approximate the $L^2(\mu)$ norm of the error as a function of the initial state of the underlying system. The iterative errors are seen to increase and eventually settle around $\sqrt{2}$, while the error from the direct mode settles at $1.0$. We show that this is a universal behavior, based on a mathematical framework \eqref{eqn:paradigm:1} that unifies both these paradigms, and both these modes. Based on this framework, we develop Theorems {\rm \ref{thm:direct}} and {\rm \ref{thm:iterative}} which provide expressions for the asymptotic behavior of these errors, and are consistent with these graphs. Also see Section~{\rm \ref{sec:conclus}} for an extended analysis. }
		\label{fig:compare}
	\end{figure}
	
	\begin{figure}[!ht]\center
		\includegraphics[width=.48\linewidth]{\figs 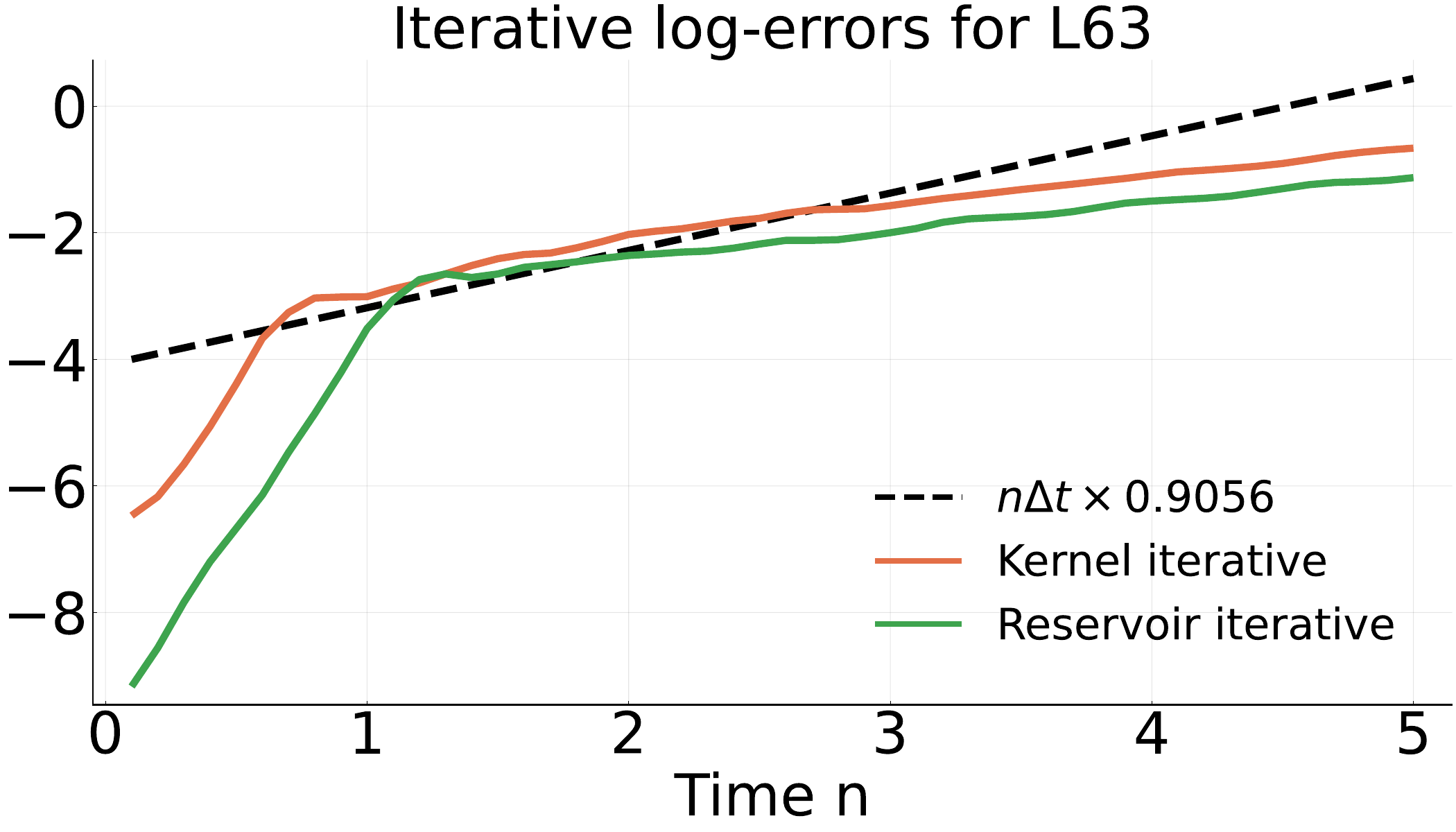}
		\includegraphics[width=.48\linewidth]{\figs 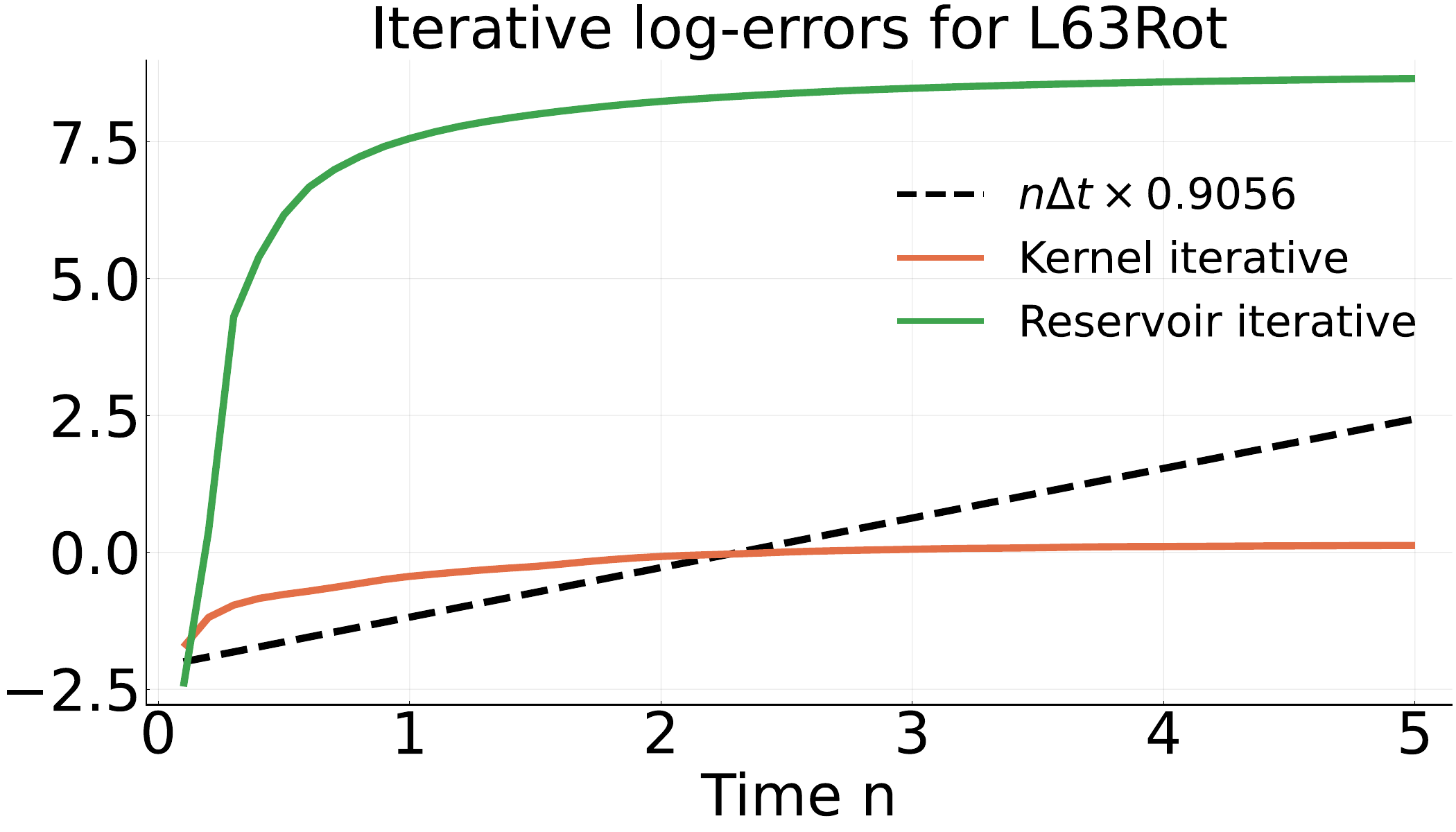}
		
		\includegraphics[width=.48\linewidth]{\figs 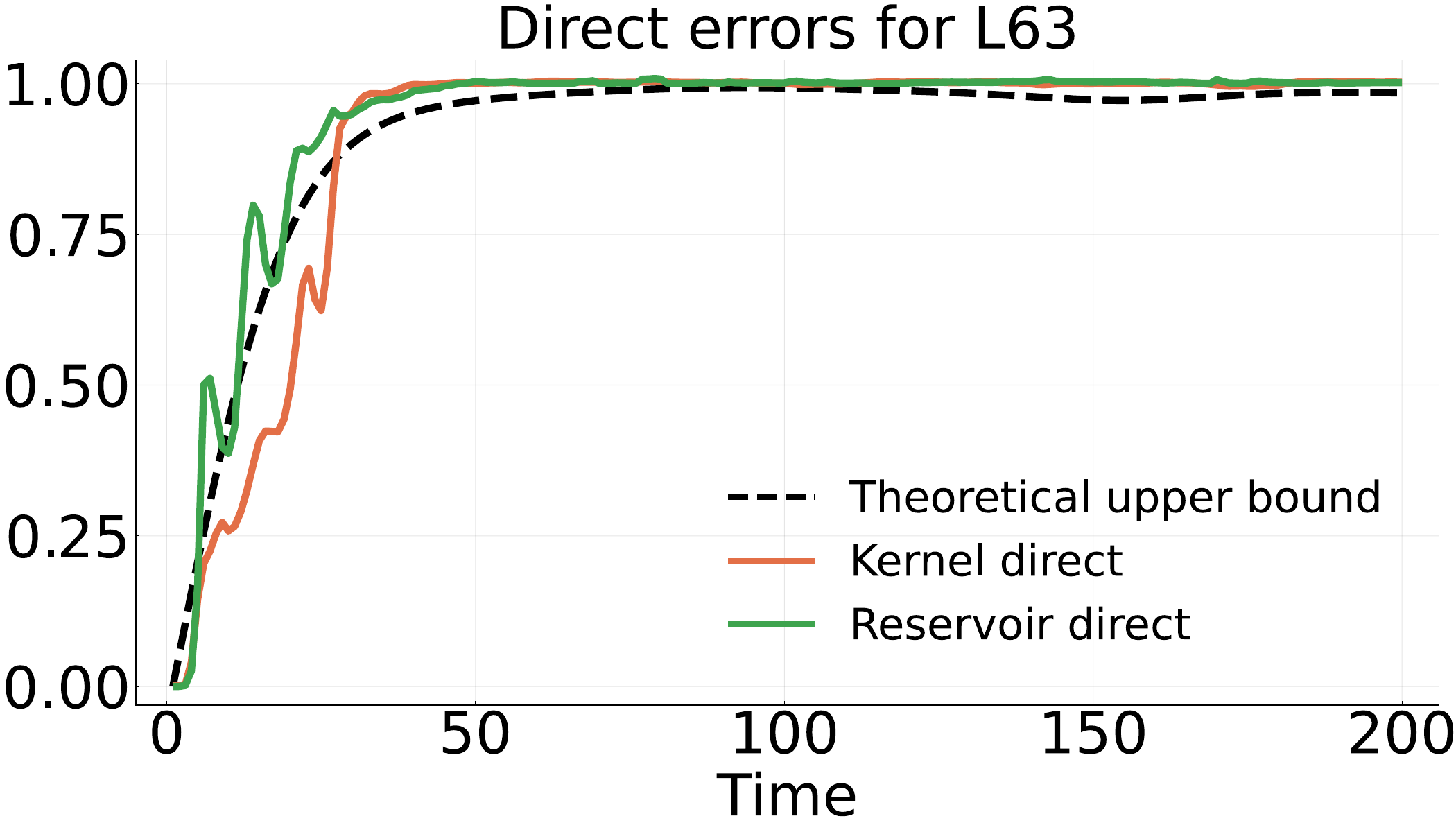}
		\includegraphics[width=.48\linewidth]{\figs 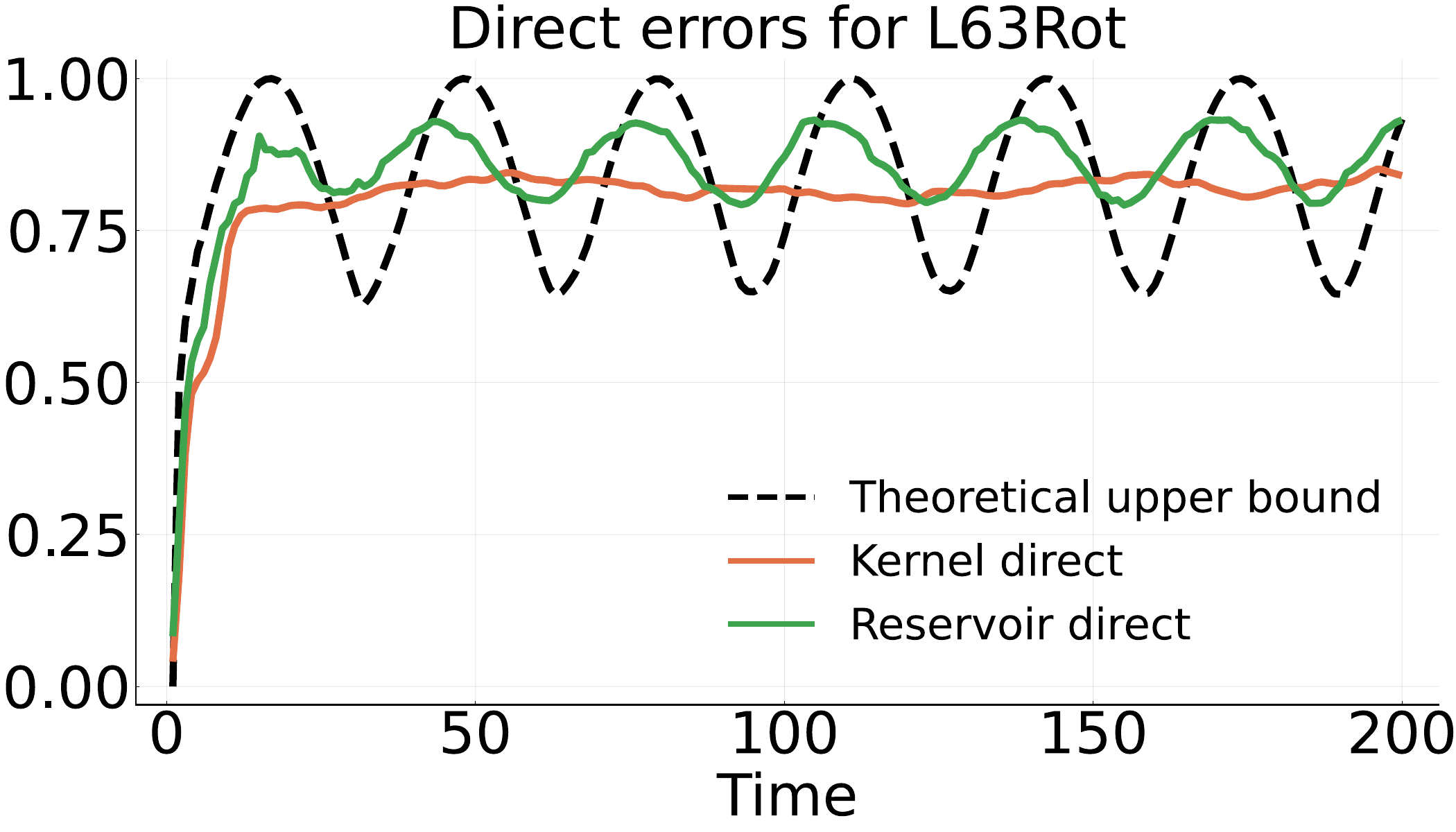}
		\caption{ Error analysis using theoretical results. The top row show how the iterative errors of the L63 and L63Rot systems compare with the theoretical bounds of Theorem~\ref{thm:iterative}. The dashed lines each have slopes $\lambda_1 \Delta t$, with $\Delta t$ being the sampling interval, and $\lambda_1 \approx 0.9056$ for both the systems. The bottom row compares the errors from direct forecast with the autocorrelation bound of \eqref{eqn:phi_hypo_autocorr}. The assumption there that $w$ lies in the hypothesis space is not met, and thus we see some fluctuations above the theoretical upper bound. Together, these plots indicate conformity with our theoretical predictions. }
		\label{fig:error_analysis}
	\end{figure}
	
	We now apply the two learning paradigms -- delay-coordinates and invariant-graphs, on three test-problems :
	\begin{enumerate}[(i)]
		\item A quasiperiodic rotation on a two-dimensional torus [Figure~\ref{fig:T2_L63Rot}, top panel].
		\begin{equation} \label{eqn:T2}
			\left( \theta^{(1)}(n+1), \theta^{(2)}(n+1) \right) = \left( \theta^{(1)}(n+1), \theta^{(2)}(n+1) \right) + ( 
			\rho_1, \rho_2 ) \bmod 2\pi .
		\end{equation}
		Here $\theta^{(1)}$ and $\theta^{(2)}$ are angular coordinates on the torus, and $( \rho_1, \rho_2 ) $ is the rotation vector. This system is analyzed in the upper panel of Figure~\ref{fig:T2_L63Rot}.
		\item The Lorenz-63 (L63) system. Let $\Phi^{t}_{\rm L63}$ denote the flow under the Lorenz63 system. Fix a sampling interval  $\Delta t$. This leads to the discrete time system :
		\begin{equation} \label{eqn:L63}
			\left( x_{n+1}, y_{n+1}, z_{n+1} \right) = \Psi^{\Delta t}_{\rm L63} \left( x_{n}, y_{n}, z_{n} \right) .
		\end{equation}
		$\Phi^{t}_{\rm L63}$ has a unique physical measure which has been proved to be nonuniformly hyperbolic and mixing. This system is analyzed in the bottom panel of Figure~\ref{fig:compare}.
		\item A dynamical system formed by taking the Cartesian product of L63 with a simple harmonic oscillator. Such a system will have a mixed spectrum, with the space $\Disc$ generated by a single base eigenfunction. We shall refer to this system as L63Rot.
		\begin{equation} \label{eqn:L63Rot} \begin{split}
				\theta_{n+1} &= \theta_n + \rho \bmod 2\pi \\
				\left( x_{n+1}, y_{n+1}, z_{n+1} \right) &= \Psi^{\Delta t}_{\rm L63} \left( x_{n}, y_{n}, z_{n} \right) .
		\end{split} \end{equation}
		This system is analyzed in the bottom panel of Figure~\ref{fig:T2_L63Rot}.
	\end{enumerate}
	Both of the two paradigms of invariant-graphs and delay-coordinate embedding were tested for these three dynamical systems. The former was implemented using reservoir systems, and the latter using kernel regression. The results of our computations in Figures \ref{fig:compare}, \ref{fig:T2_L63Rot}, and \ref{fig:error_analysis} illustrate the consequences of Theorem~\ref{thm:direct}, Theorem~\ref{thm:iterative}. 
	
	Figure~\ref{fig:error_analysis} highlights the two most important conclusions from our results. Firstly, as seen in the top row, the iterative errors grow at an exponential rate comparable to the top Lyapunov exponent $\lambda_1$. Secondly, if the hypothesis space is large enough, then the direct error is bounded above by the formula \eqref{eqn:phi_hypo_autocorr}. The errors from the iterative forecasts made using the reservoir blow up. This is because the standard reservoir computers are not guaranteed to be stable. 
	
	\paragraph{Asymptotic analysis} Theorem~\ref{thm:iterative} gives an upper bound for the long term behavior of the iterative error. The theoretical bound for exponential rate of growth is indicated by the slope of the black dashed line, and is $\approx 0.9056\Delta t$. So although the initial exponential rate of errors seem to be larger than this, by choosing a multiplicative constant large enough, the error graph still remains underneath the theoretical curve. The offset of the straight dashed line equals the logarithm of this multiplicative constant. Thus as long as the long-term averaged error growth rate is less than $\approx 0.9056\Delta t$, there will always be a multiplicative constant large enough to satisfy the bounds in  \eqref{eqn:def:iter_bound_L2}.
	
	\paragraph{Ideal direct error} The errors from the direct error occasionally cross the theoretical bound indicated by the black dashed line. This is because, the bound in \eqref{eqn:phi_hypo_autocorr} assumes that the learning error for $w$ is zero, i.e., $w$ lies in the hypothesis space $\mathcal{H}$. In most situations such as in our experiments, there is always a small learning error. An extended analysis for this situation is an interesting and open task.
	
	This completes the experimental verification of our theoretical results.
	
	\section{Conclusions} \label{sec:conclus}
	
	The learning theory of dynamical systems poses many challenges not usually encountered in other areas of Approximation theory. As discussed in Section \ref{sec:intro} although a dynamical system is a map $f$, it is overall an iterative procedure. So most of the salient features of a dynamical system are asymptotic, such as an attractor and an invariant measure. Other features are induced in associated spaces, such as the Koopman or transfer operator. The dependence of these features on the dynamics map is not necessarily continuous. As a result approximation of the dynamics map alone does not necessarily guarantee an approximation of the other features. The present article takes a dive into these challenges of learning theory in the context of dynamical systems.
	
	The vast field of data-driven dynamical systems thus needs to be separated based on the objective being targeted. The article presents four such objectives. The diagram in Figure \ref{fig:outline1} is broad enough to contain most learning techniques. These techniques are spread across different branches of applied Dynamical systems such as medical imaging, signal processing, and climate science. The challenge faced by most practitioners of these techniques is due to the fact that a numerical method built to approximate one aspect of the dynamics may not be adequate to approximate a different aspect. The reason for this comes back to the various dis-connections between these aspects. For example, prediction techniques are limited because they employ an estimate of the dynamics law to estimate the Koopman operator. Similarly, A purely topological approach to the dataset fails to capture the intricacies of the dynamics. Figures \ref{fig:outline2} and \ref{fig:outline3} present some examples when these cross connections are feasible, and when they are not. 
	
	We focused on the basic task of learning $f$ in more depth in Section \ref{sec:learn}. There has also been a huge body of work focused solely on the task of reconstructing $f$ from data. Table \ref{tab:learn} provides a broad survey of many of these techniques. An important realization of Section \ref{sec:learn} is that all of these fall under an abstract framework presented jointly by Assumption \ref{A:data} and \ref{A:pPhi}. The ingredient of this framework is a map $g$ which we have called the \emph{embedding mechanism}. The two main embedding mechanisms are invariant-graphs and delay-coordinates, and these methods are compared in Table \ref{tab:prdgm}.
	
	We have examined in detail the difficulty of estimating the invariant set from an an approximation $\hat{f}$ of $f$; and the difficulty of estimating the Koopman operator from $\hat{f}$. The former is examined in Section \ref{sec:stability}. It occurs due to the creation of additional Lyapunov exponents during a data-driven reconstruction. The latter is examined in Section \ref{sec:predict}.
	
	The top Lyapunov exponent $\lambda_1(\hat{T})$ of the reconstructed model $\hat{T}$ depends on the behavior of the feedback function in a neighborhood of the image of the attractor. This exponent characterizes the stability of the reconstructed invariant set. Theorems \ref{thm:lambda1} and \ref{thm:iterative} suggest that instead of measuring the overfitting error via the usual smoothness bound, a good candidate would be to take into account the behavior of $\bar{w}$ in a neighborhood of the dataset. 
	
	We have examined in depth the task of forecasting, using a reconstruction of the dynamics law. It is important to realize that forecasting can be done in at least two ways - direct and iterative. Operator theoretically, the direct method is essentially an application of $\pi U^n \phi$ and the iterative is $U^{n-1}\pi U\phi$.  The fundamental difference in these forecasting methods is the location of the projection, $\pi$, in the order of operations being estimated. Which method achieves better results depends on the sensitivity and mixing characteristics of the dynamics.
	
	Equation \eqref{eqn:direct:1} is an operator theoretic restatement of the formula \eqref{eqn:def:err_direct}. Theorem~\ref{thm:direct} relates the growth of $\text{error}_{\text{direct}}$ with the rate of decay of correlations. The performance of the direct forecast is tied to the spectral property of the Koopman operator. On the other hand the performance of the iterative forecast is tied to the Lyapunov spectrum of the dynamics as shown in Theorem~\ref{thm:iterative}. This bivalence is displayed by all dynamical systems.
	
	Thus the task of data-driven discovery of dynamical systems is multi-faceted, and its limits are tied to properties inherent to the dynamics, not just the availability of data and the power of numerical methods. A proper understanding of these limitations is important to build and assess numerical methods.
	
	\bibliographystyle{unsrt}
	\bibliography{References,SIADS_bib}
\end{document}